\documentclass{article}
\usepackage[T1]{fontenc}
\usepackage[utf8]{inputenc}
\usepackage{amsmath,amssymb,amsthm,enumitem,graphicx,float,subcaption,titlesec,xcolor,tikz,mathrsfs,mathtools}
\usepackage[left=0.75 in, right=0.75 in, top = 0.75 in, bottom = 0.75 in]{geometry}
\usepackage[justification=centering]{caption}
\usepackage[english]{babel}
\usepackage{pifont,abstract}
\usepackage[hidelinks,colorlinks=true,linkcolor=black,citecolor=black,filecolor=black,urlcolor=cyan]{hyperref}
\usetikzlibrary{decorations.markings,matrix}

\def\papertitle{Construction and Applications of Trisections of Low Genus on Del Pezzo Surfaces of Degree One}
\def\theauthor{Julie Desjardins and Vojin Jovanovic}

\definecolor{highlighted}{RGB}{255,0,0}
\definecolor{functioncolor}{RGB}{192,63,255}
\definecolor{fibrecolor}{RGB}{25,155,250}

\renewcommand\maketitle{%
    \begin{center}%
    \noindent\vphantom{0pt}%
    \vspace{50pt}%
    
    {%
    \noindent\bfseries\Large\papertitle{}%
    }\vspace{15pt}%
    
    {%
    \noindent\large\theauthor{}%
    }\vspace{2ex}%
    \end{center}%
}

\theoremstyle{definition}
\newtheorem{definition}{Definition}[section]
\newtheorem{proposition}[definition]{Proposition}
\newtheorem{corollary}[definition]{Corollary}
\newtheorem{lemma}[definition]{Lemma}
\newtheorem{theorem}[definition]{Theorem}
\newtheorem*{key}{Key Concept}

\theoremstyle{remark}
\newtheorem*{remark}{Remark}
\newtheorem{example}[definition]{Example}

\newcommand{\nin}{\not\in}

\newcommand{\Q}{\mathbb{Q}}
\newcommand{\Z}{\mathbb{Z}}
\newcommand{\N}{\mathbb{N}}
\newcommand{\E}{\mathcal{E}}
\newcommand{\F}{\mathcal{F}}
\newcommand{\T}{\mathcal{T}}
\renewcommand{\P}{\mathbb{P}}

\def\xscale{0.1}
\def\yscale{0.4}
\def\xoff{0.241526}
\def\sample{25}

\def\del{\partial}
\def\magma{\texttt{magma}}
\def\famL{\text{\textbf{Fam}}_L}
\def\famLk{\text{\textbf{Fam}}_k}
\def\famk{\text{\textbf{fam}}_k}
\def\famtwo{\text{\textbf{Fam}}_2}
\def\famone{\textbf{Fam}_1}

\def\leftoffsetepsilon{5 em}
\def\rulefill{\hspace{\textwidth}}

\def\codelink{https://drive.google.com/drive/folders/1eA5K63udWjmSaPtLolZpMH-4pJ0cI8oK?usp=drive_link}
\def\thecode{\href{\codelink}{\underline{here}}; alternatively, the url can be found at the bottom of the paper}

\begin{document}

\maketitle

\begin{abstract}
    Consider a rational elliptic surface over a field $k$ with characteristic $0$ given by $\E: y^2 = x^3 + f(t)x + g(t)$, with $f,g\in k[t]$, $\deg(f) \leq 4$ and $\deg(g) \leq 6$. 
    If all the bad fibres are irreducible, such a surface comes from the blow-up of a del Pezzo surface of degree one.
    We are interested in studying multisections, curves which intersect each fibre a fixed number of times, specifically, trisections (three times).
    Many configurations of singularities on a trisection lead to a lower genus.
    Here, we focus on of several them: by specifying conditions on the coefficients $f,g$ of the surface $\E$, and looking at trisections which pass through a given point three times, we obtain a pencil of cubics on such surfaces.
    Our construction allows us to prove in several cases the Zariski density of the rational points. 
    This is especially interesting since the results in this regard are partial for del Pezzo surfaces of degree one.
\end{abstract}
\section{Introduction}
\label{sec:intro}

A \textit{del Pezzo surface} $S$ is a smooth, projective, and geometrically integral surface defined over a field $k$, with ample anti-canonical divisor, $-K_S$.
The \textit{degree of a del Pezzo surface} is the self-intersection number of its canonical divisor, which is an integer $d$ between 1 and 9.
Every del Pezzo surface is isomorphic to the blow-up of $\P^2$ at $9-d$ points in general position if $d\neq 8$, or isomorphic to $\P\times \P$ when $d = 8$ \cite[Theorem 24.4]{Manin_1986}.
A \textit{rational elliptic surface} $\E$ is a projective surface with a rational elliptic fibration, i.e., a rational proper morphism $\pi : \E \rightarrow \P^1$ with connected fibres, whose generic fibre is a smooth curve of genus one.
Iskovskih shows that over an algebraically closed field, each rational elliptic surface has a minimal model which is a del Pezzo surface or a conic bundle \cite{Iskovskih_1980}.

A variety $X$ is said to be \textit{$k$-unirational} if there exists a dominant rational map $\P_k^n \rightarrow X$.
For an arbitrary field $k$, a del Pezzo surface of degree $d\geq 2$ is $k$-unirational if it contains a $k$-rational point, with additional technical conditions in the case of $d = 2$ and $k$ an infinite field.
Segre shows this for degree 3 and $k = \Q$ in \cite{Segre_1943} and \cite{Segre_1951}.
Manin studies the complete case for degree $\geq 5$, and provides partial results for degrees $3$ and $4$ for large-enough cardinality of $k$ in \cite{Manin_1986}.
Kollár shows the remaining case of $d = 3$ in \cite{Kollár_2002}; see \cite[Proposition 5.19]{Pieropan_2012} for the general case of $d = 4$.
In \cite{STVA_2014}, the authors study the case of arbitrary fields $k$ with characteristic other than 2 and del Pezzo surfaces of degree $d = 2$, with the additional assumption that the $k$-rational point is not contained in the ramification divisor, nor in four exceptional curves.
The same paper shows that all del Pezzo surfaces of degree two are $k$-unirational when $k$ is an arbitrary finite field, except for the cases of three explicit surfaces, which are covered by \cite{FvL_2015}.

Birational equivalence preserves $k$-unirationality, so any surface which is birationally equivalent to a del Pezzo surface of degree $\geq 2$ as above is $k$-unirational.
A del Pezzo surface $S$ has \textit{minimal degree} $d$ if it is not birationally equivalent to any del Pezzo surface of degree $> d$.
In \cite{KM_2017}, the authors show that del Pezzo surfaces of degree one, which arise from a conic bundle, are $k$-unirational.
Such surfaces all have Picard rank two, however, the question of $k$-unirationality remains open with respect to minimal del Pezzo surfaces of degree one and Picard rank one.
To the authors' knowledge, there are no examples of such del Pezzo surfaces which are either known to be $k$-unirational or known not to be $k$-unirational.

For infinite fields $k$, $k$-unirationality of a variety $X$ implies that the $k$-rational points $X(k)$ are dense in the Zariski topology.
So, it is natural to study this weaker property on del Pezzo surfaces of degree one, with particular interest in the case where the Picard rank is one.
Any del Pezzo surface of degree one can be blown-up at the anti-canonical base point, giving rise to an elliptic fibration on the resulting surface.
This rational elliptic surface has only irreducible bad fibres; conversely, such a rational elliptic surface comes from a blow-up of a del Pezzo surface of degree one \cite{Desjardins_2016}.
This is discussed in more detail in section \ref{sec:prerequisites}.
Regarding del Pezzo surfaces of degree one, several results are known.
In \cite{Ulas_2008}, Ulas constructs a 12-section with a genus one component on the family of del Pezzo surfaces of degree one which correspond to the following elliptic surfaces:
\begin{equation*}
    S_f : y^2 = x^3 + f(t), 
\end{equation*}
where $f(t) = t^5 + at^3 + bt^2 + ct + d\in\Z[t]$.
In \cite{SvL_2014}, the authors extend this construction to an arbitrary del Pezzo surface of degree one over an infinite field $k$ of characteristic different from 2 or 3. 
If further technical criteria are met on the surface, this section is used to conclude that the surface has Zariski dense rational points.
In \cite{Bulthuis_2018}, Bulthuis, in a master's thesis supervised by van Luijk, shows that any del Pezzo surface whose corresponding rational elliptic surface contains a smooth fibre with a torsion point also contains a pencil of genus one curves.
In \cite{DW_2022} the authors construct a pencil of trisections, all of which have genus one or lower, on the family of del Pezzo surfaces of degree one given by:
\begin{equation*}
    Y^2 = X^3 + aZ^6 + bZ^3W^3 + cW^6,
\end{equation*}
defined over a field $k$ of characteristic zero, where $a,c$ are both non-zero.
These trisections are constructed by specifying conditions on the coefficients to force a total of three singularities, thus lowering the genus.
This family is then used to show that the surfaces all have Zariski dense rational points.
A similar method is used in \cite{Nijgh_2022} for the family of del Pezzo surfaces of degree one given by:
\begin{equation*}
    Y^2 + c_1XYW + \left[c_2f(Z/W) + c_3\right]YW^3 = X^3 + c_4X^2W^2 + \left[c_5f(Z/W) + c_6\right]XW^4 + \left[c_7f(Z/W)^2 + c_8f(Z/W) + c_9\right] W^6, 
\end{equation*}
where $f\in k[t]$ is a degree three polynomial over a field of characteristic zero.
These papers are discussed in more detail in section \ref{sec:known-results}.

In this paper, we construct trisections of low genus.
For the rest of the paper let $k$ be an arbitrary field of characteristic zero.

\begin{theorem}
    \label{thm:intro-main}
    Let $S/k$ be a del Pezzo surface of degree one given by the form:
    \begin{equation*}
        S : Y^2 = X^3 + F(Z,W)X + G(Z,W)
    \end{equation*}
    in the weighted projective space $\P(2,3,1,1)$, where $F,G\in k[Z,W]$ are homogeneous polynomials of degrees four and six, respectively.
    Let $f(t) := F(t,1)$, and $g(t) := G(t,1)$, and suppose there exists a field extension $L/k$, along with $a,b,t_0\in L$, such that the following all hold:
    \begin{enumerate}[label=(\roman*)]
        \item $t_0(at_0 + b) \neq 0$; 

        \item $f(t_0) \neq -(at_0 + b)^4/3$, and;

        \item the polynomial
        \begin{equation*}
            P(t) = -f(t)^2/4 + (at+b)^4f(t)/6 + (at+b)^2g(t) + (at+b)^8/108
        \end{equation*}
        has a double root at $t = t_0$.
    \end{enumerate}
    Then, the $L$-rational point $Q = ((at_0 + b)^2/3 : (at_0 + b)^3/6 + f(t_0)/(2(at_0 + b)) : t_0 : 1)$ lies on the surface $S$, and there exists a trisection on $S$, defined over $L$, with a triple singularity at $Q$.
    
    Assume there exists a fibre $\F$ with $Q\nin \F$, and $|\F(L)| = \infty$; then, there exists a pencil of trisections on $S$, defined over $L$, which all have a triple singularity at $Q$.
    Moreover, each of these trisections has genus at most one.
    The converse may not hold.
\end{theorem}

This pencil of trisections can have several uses, however, it is useful to consider surfaces where this pencil is defined over the base field instead of an extension.
Some surfaces with this property are provided by the following corollary.

\begin{corollary}
    Let $S/k$ be defined as above and suppose that there exists an $a\in k^*$ such that the polynomials
    \begin{align*}
        \hat f(t) := f(t) - a^4t^4/3 & & \hat g(t) := g(t) + a^6t^6/27
    \end{align*}
    have a shared root $t_0\in k^*$, with $t_0$ being a double root of $\hat g$.
    Then, the trisections from Theorem \ref{thm:intro-main} are defined over $k$.
\end{corollary}

\begin{example}
    The following family of rational elliptic surfaces $\E_\alpha$ is generically the blow-up of del Pezzo surfaces of degree one which satisfy the conditions of Theorem \ref{thm:intro-main}, with $L = k$ and $a = 1, b = 0, t_0 = 1$:
    \begin{equation*}
        \E_\alpha : 
        \begin{cases}
            f(t) = \frac{4}{3}t^4 - 1 \\
            g(t) = \left(\alpha - \frac{1}{27}\right)t^6 + \left(-2\alpha + 1\right)t^5 + \left(\alpha - 1\right)t^4 + \left(-16\alpha - \frac{890}{27}\right)t^2 + \left(32\alpha + \frac{1753}{27}\right)t - 16\alpha - \frac{863}{27}.
        \end{cases}
    \end{equation*}
    The pencil of trisections is given by:
    \begin{equation*}
        T_h : xt + \left(-\frac{45}{2}\alpha - h - \frac{454}{9}\right)t^3 + \left(45\alpha + 3h + \frac{926}{9}\right)t^2 + \left(-\frac{45}{2}\alpha - 3h - \frac{472}{9}\right)t + h \subset \E_\alpha,
    \end{equation*}
    each of which has a triple singularity at the point $Q = (1/3, 1/3, 1)$.
    If we set $h = -45\alpha - 863/9$, we get the point $R = (1, 1, 2)\in T_h$.
    This example and others are shown in section \ref{sec:examples}.
\end{example}

We can use these trisections to prove the Zariski density of the $L$-rational points of $S$ in the following case.

\begin{theorem}
    \label{thm:intro-ZD}
    Let $S$ be as in Theorem \ref{thm:intro-main}, and suppose there is a fibre $\F\subset S$, with $Q\nin \F$, such that $|\F(L)| = \infty$.
    Further, suppose either:
    \begin{enumerate}[label=(\roman*)]
        \item one of the trisections contains a singularity at a point other than $Q$, or;
    
        \item infinitely many of the genus one curves in the pencil contain an $L$-rational point of infinite order, or;

        \item there exists a genus one curve in the pencil which contains infinitely many $L$-rational points which are non-torsion on their respective fibres.
    \end{enumerate}
    Then, the $L$-rational points, $S(L)$, are Zariski dense in $S$.
\end{theorem}

Corresponding to the first case, we have the following theorem.

\begin{theorem}
    \label{thm:intro-DP1}
    There exists an infinite family of rational elliptic surfaces $\E$ defined over $k$, with corresponding minimal model $S$.
    Each $\E$ in this family contains a trisection with a triple singularity at a point $Q$, and a double singularity at a point $R \neq Q$.
    Each of these surfaces has Zariski dense $k$-rational points.
    Generically this family only has irreducible bad fibres.
\end{theorem}

\begin{remark}
    The minimal models $S$ of the above family are all members of the family defined in Theorem \ref{thm:intro-main}.
\end{remark}

This leads to the natural question of whether one can construct a pencil of trisections, all of which have an additional singularity.
If this holds, we are left with either a conic or the surface will have a one dimensional family of genus zero curves, allowing us to show unirationality.
We believe this is not possible without increasing the minimal degree of the del Pezzo surface.
\cite[Proposition 2.9]{Gao_2023} shows that this is not possible in the case of sections or bisections.
While we cannot use this method to prove unirationality for a del Pezzo surface of minimal degree one, we can apply it to prove the following theorem about a family of del Pezzo surfaces, all of which have minimal degree of at least two.

\begin{theorem}
    \label{thm:intro-DP2}
    Let $\E$ be the rational elliptic surface given by:
    \begin{equation*}
        \E : y^2 = x^3 + f(t)x + g(t),
    \end{equation*}
    with $f(t) = \sum_{i=0}^4 f_it^i$, and $g(t) = \sum_{i=0}^6 g_it^i$.
    Further suppose that $f_0,f_1,g_0,g_1,g_2,g_3,g_4$ all satisfy equations (\ref{eq:constraint-F0} - \ref{eq:constraint-G4}).
    Then, $\E$ has a minimal model $S$, which is a del Pezzo surface of degree at least two; $S$ is unirational and, further, the degree of the rational map $\P^2 \rightarrow S$ is three.
\end{theorem}

The interest of this result, compared to \cite{STVA_2014}, is we obtain a low degree rational map.
The family in Theorem \ref{thm:intro-DP2} has ten free variables, each element of the family has a minimal model which is a del Pezzo surface.
The following is an example of one of the surfaces in this family.

\begin{example}
    The rational elliptic surface $\E$ with coefficients given by
    \begin{equation*}
        \E :
        \begin{cases}
            f(t) = t^4 + t^3 + t^2 - \frac{68}{3}t + \frac{43}{3}, \\
            g(t) = t^6 + t^5 - \frac{2783}{144}t^4 + \frac{4175}{216}t^3 + \frac{501}{16}t^2 - \frac{433}{12}t + \frac{811}{108}
        \end{cases}
    \end{equation*}
    contains a pencil of genus zero trisections defined as:
    \begin{equation*}
        T_h : y = xt + x + \left(\frac{13}{4} - \frac{1}{2}h\right)t^3 + \left(2h - \frac{103}{12}\right)t^2 + \left(\frac{8}{3} - \frac{5}{2}h\right)t + h.
    \end{equation*}
    Each trisection has a triple singularity at the point $Q = (4/3, 0, 1)$, and a double singularity at $R = (1, 0, 2)$.
\end{example}
\section{Prerequisites}
\label{sec:prerequisites}

\begin{definition}
    A \textit{del Pezzo surface} $S$ is a smooth, projective, and geometrically integral surface over $k$ with ample anti-canonical divisor $-K_S$.
    The \textit{degree of a del Pezzo surface} is the self intersection number of the canonical divisor $K_S$, which is an integer between 1 and 9.
\end{definition}

\begin{definition}
    A surface $X$ over $k$ is said to be \textit{$k$-unirational} if there exists a dominant map $\P^n_k \rightarrow X$, for some natural number $n$.
\end{definition}

For the purposes of this paper, we instead use an equivalent definition of a del Pezzo surface of degree one (DP1); see \cite[Theorem 3.3.5]{Kollar_1996} for details.

\begin{definition}
    \label{def:proj-surface}
    Let $k$ be a field of characteristic zero.
    Consider the surface in the weighted projective space $\P(2,3,1,1)$ with corresponding variables $(X:Y:Z:W)$ given by:
    \begin{equation*}
        S : Y^2 = X^3 + F(Z,W)X + G(Z,W),
    \end{equation*}
    where $F,G\in k[Z,W]$ are homogeneous polynomials in two variables, with $\deg(F) = 4$ and $\deg(G) = 6$.
    If this surface is smooth, then it is a del Pezzo surface of degree one.
    After a blow up at $\infty$, we are left with a rational elliptic surface in affine space with corresponding variables $(x,y,t)$ of the form:
    \begin{equation}
        \label{eq:affinesurface}
        \E : y^2 = x^3 + f(t)x + g(t),
    \end{equation}
    where $t = Z/W$, $f(t) = F(t,1)$, and $g(t) = G(t,1)$.
    This rational elliptic surface has only irreducible bad fibres.
    There is, in fact, a correspondence between these two concepts: if a rational elliptic surface of the form found in equation (\ref{eq:affinesurface}) has only irreducible bad fibres, then it is the result of a blow-up of a del Pezzo surface of degree one \cite{Desjardins_2016}.
    This correspondence will be mentioned throughout the paper.
\end{definition}

\begin{definition}
    \label{def:multisection}
    On the surface $\E$, a \textit{multisection of degree} $n$, called an \textit{$n$-section}, is a divisor in the linear system $|-nK_S|$, where $K_S$ is the canonical divisor of $S$.
    Each such divisor is associated with an algebraic set of dimension one, hence each of the irreducible components are therefore (potentially singular) curves on the surface.
    We are interested in $3$-sections, also called trisections, i.e., multisections of degree 3.
    A general trisection is the subvariety of $\E$ corresponding to an equation of the following form:
    \begin{equation}
        \label{eq:trisection}
        C : y = axt + bx + ct^3 + dt^2 + et + h,
    \end{equation}
    where $a,b,c,d,e,h\in k$ are constants.
\end{definition}

We can obtain the arithmetic genus of a general curve $C$ using \cite[Theorem 9.1.37]{LE_2002}, which gives the following formula:

\begin{equation}
    \label{eq:Liu_genus}
    g_a(C) = 1 + \frac{C\cdot(C + K_S)}{2},
\end{equation}
where "$\cdot$" denotes the intersection number and $K_S$ is the canonical divisor of $S$.

\begin{theorem}
    \label{thm:sing-genus}
    Let $C$ be a singular curve, and let $\tilde C$ be the desingularized curve.
    For any singularity $R$ on $C$, let $r$ be the multiplicity of the singularity.
    Then, we obtain the following formula:
    \begin{equation*}
        g_a(C) = g_a(\tilde C) - \smashoperator{\sum_{\substack{R\;\in\;\text{sing}(C)}}}\frac{r(r-1)}{2}.
    \end{equation*}
\end{theorem}

\proof Apply \cite[Theorem V.3.7]{Hartshorne_1977} repeatedly.
$\qed$

In particular, for a non-singular trisection on a del Pezzo surface $S$ of degree one, equation (\ref{eq:Liu_genus}) can be calculated to be 4.
So, for a singular trisection $C$, the genus formula becomes:
\begin{equation}
    \label{eq:genus_formula}
    g_a(C) = 4 - \smashoperator{\sum_{\substack{R\;\in\;\text{sing}(C)}}}\frac{r(r-1)}{2}.
\end{equation}
In this paper, we construct families of trisections wherein each member has a triple singularity.
By (\ref{eq:genus_formula}), this singularity reduces the genus to one or lower.

\begin{key}
    Adding a singularity to a multisection will lower the genus.
    This method is used to obtain a pencil of genus one or zero curves, which can then be used to show Zariski density of the $k$-rational points on some families of surfaces.
\end{key}
\subsection{Select Known Results}
\label{sec:known-results}

There is a connection between the geometric question of the existence of a low genus multisection and the algebraic question relating to the Zariski density of the rational points on a surface.
There are many examples from prior known results where multisection components of low genus are used in part to show the Zariski density of the rational points.
This section discusses select examples to demonstrate this principle and motivate the further study of this topic.
In particular, we summarize the usage of multisections of degree twelve and degree three in past literature.

In \cite{Ulas_2008}, Ulas looks at surfaces of the form $x^2 = y^3 + f(z)$, where $f(z) = z^5 + az^3 + bz^2 + cz + d\in\Z[z]$.
Ulas shows that if the curve $Y^2 = X^3 + 135(2a - 15)X - 1350(5a + 2b - 26)$ contains infinitely many rational points, then the rational points on the original surface are Zariski dense.
Ulas constructs a 12-section on the original surface, which has a component that has genus one.
This component is associated with the above curve; in particular, a rational point on the curve corresponds to a rational point on the multisection.
If the curve has infinitely many rational points, then the multisection component has positive rank, so the rational points on the original surface are Zariski dense.

In \cite{SvL_2014}, the authors consider an arbitrary del Pezzo surface of degree one, given as in Definition \ref{def:proj-surface}, except over an infinite field $k$ of characteristic not equal to 2 or 3.
The authors construct a 12-section on this surface and use this multisection to conclude Zariski density of the rational points so long as several additional properties also hold.

In \cite{DW_2022}, the authors consider the surface $Y^2 = X^3 + aZ^6 + bZ^3W^3 + cW^6$ defined over a field $k$ of characteristic zero, where $a,c$ are both non-zero.
The authors show that if $S$ contains a rational point with non-zero $Z,W$-coordinates such that the corresponding point on $\E$ is torsion-free on its fiber, then the $k$-rational points are Zariski dense in $S$.
Here, $\E$ is the elliptic surface obtained by blowing up the base point of the anti-canonical linear system $|-K_S|$.
The proof proceeds by constructing a family of trisections on $\E$ and forcing each of them to contain three singularities, each of multiplicity two.
This leads to the following three cases:
\begin{enumerate}[label=(\roman*)]
    \item a member of the family contains a section defined over $k$, or;

    \item a member of the family is geometrically integral of genus zero, or;

    \item infinitely many members are geometrically integral of genus one. 
\end{enumerate}
In each of the three cases the authors show that $S(k)$ is Zariski dense in $S$.

In \cite{Nijgh_2022}, Nijgh, in a master's thesis supervised by van Luijk, looks at smooth surfaces of the form:
\begin{equation*}
    Y^2 + c_1XYW + \left[c_2f(Z/W) + c_3\right]YW^3 = X^3 + c_4X^2W^2 + \left[c_5f(Z/W) + c_6\right]XW^4 + \left[c_7f(Z/W)^2 + c_8f(Z/W) + c_9\right] W^6, 
\end{equation*}
where $f(t) = f_3t^3 + f_2t^2 + f_1t + f_0\in k[t]$ is a degree three polynomial over a field of characteristic zero.
Nijgh shows that if there exists an element $t_0\in k$, such that: (i) $f-f(t_0)$ is separable, (ii) $3t_0 \neq -f_2/f_3$, and (iii) for $\F := Z - t_0W = 0\subset S$, $|\F(k)| = \infty$, then the $k$-rational points in $S$ are Zariski dense.
When $k$ is of finite type over $\Q$, the converse holds as well.
The proof proceeds by constructing a pencil of trisections, each with at least three singularities of multiplicity two, meaning each has genus at most one.
Nijgh then constructs related surfaces to this pencil, and using a rational dominant map into $S$, shows that the $k$-rational points of $S$ are Zariski dense by looking at the Zariski density of the $k$-rational points on the pencil and related surfaces.
\section{Constructing the Pencil}
\label{sec:main}

We start by defining two families of DP1s as follows.

\begin{definition}
    \label{def:fam-L}
    For a fixed field $k$ of characteristic zero, fix a field extension $L/k$, and define the following family.
    Let $\famL$ be a family of all del Pezzo surfaces $S$ of degree one which satisfy the following condition: there exist fixed elements $a,b,t_0\in L$, with $t_0 \neq 0$, $(at_0 + b) \neq 0$ and $f(t_0) \neq -(at_0 + b)^4/3$, such that $t_0$ is a double root of the polynomial
    \begin{equation}
        \label{eq:condition}
        P(t) = -f(t)^2/4 + (at+b)^4f(t)/6 + (at+b)^2g(t) + (at+b)^8/108.
    \end{equation}
\end{definition}

\begin{remark}
    If we do not assume that $f(t_0) \neq -(at_0 + b)^4/3$, then the surface $S$ in $\famL$ will not be a del Pezzo surface of minimal degree one.
    See Lemma \ref{lemma:f-t_0} for details.
\end{remark}

\begin{definition}
    \label{def:fam-k}
    Once again for a fixed field $k$ of characteristic zero, let $\famk$ be the family of all del Pezzo surfaces $S$ of degree one which satisfy the following condition:
    There exists fixed non-zero $a,t_0\in k^*$ such that $f(t) - a^4t^4/3$ has a root at $t_0$, and $g(t) + a^6t^t/27$ has a double root at $t_0$.
\end{definition}

\begin{remark}
    We note that $\famk$, is not the same family as $\famLk$, with the latter being Definition \ref{def:fam-L} for the extension $L=k/k$.
    However, $\famk\subset\famLk\subset\famL$ for any field extension $L$.
    Fixing $b=0$, and an arbitrary field extension $L/k$, then $P(t)$ becomes:
    \begin{equation*}
        P(t) = -f(t)^2/4 + (at)^4f(t)/6 + (at)^2g(t) + (at)^8/108.
    \end{equation*}
    By completing the square with respect to $f(t)$, we get:
    \begin{equation*}
        P(t) = -\frac{\left(f(t) - (at)^4/3\right)^2}{4} + (at)^2\left(g(t) + (at)^6/27\right).
    \end{equation*}
    Let $S\in\famk$; from Definition \ref{def:fam-k}, we have $a,t_0\in k^*$ with $f(t) - a^4t^4/3$ having a root at $t_0$, and $g(t) + a^6t^6/27$ having a double root at $t_0$.
    So, $a,b,t_0\in k\subset L$, and by the above equation, we see that $P(t)$ has a double root at $t_0$.
    Also notice that $f(t_0) = a^4t_0^4/3 \neq -a^4t_0^4/3$, since $a,t_0\neq 0$, thus $S\in\famLk\subset\famL$.
\end{remark}

\begin{lemma}
    \label{lemma:trisection}
    For all field extensions $L/k$, and for all $S\in\famL$, $S$ contains a family of trisections over $L$ with a single free variable.
    Each of the trisections has a triple singularity at the point $Q = ((at_0 + b)^2/3,(at_0 + b)^3/6 + f(t_0)/(2(at_0 + b)),t_0)$.
\end{lemma}

\proof

Each surface in $\famL$ contains the family of trisections given by
\begin{equation*}
    T_\gamma : y = axt + bt + c(\gamma)t^3 + d(\gamma)t^2 + e(\gamma)t + \gamma \subset \E,
\end{equation*}
where $\gamma$ is the parameterizing variable and the variables $c(\gamma), d(\gamma)$, and $e(\gamma)$ are defined as follows:
\begin{multline}
    \label{eq:result-c}
    c(\gamma) = \frac{1}{2t_0^3(at_0 + b)^3(f(t_0) + (at_0 + b)^4/3)}\bigg(-(at_0 + b)^8\left((b - at_0/2)^2 - 7a^2t_0^2/4\right)/9 \\
    - 2\gamma(at_0 + b)^7/3 + (at_0 + b)^4\left(c_1(t_0) + t_0^2g''(t_0)\right) - 2\gamma(at_0 + b)^3f(t_0) + c_2(t_0)\bigg),
\end{multline}
where $c_1(t_0)$ and $c_2(t_0)$ are polynomials over $L$, given by the equations (\ref{eq:c1}) and (\ref{eq:c2}) below;
\begin{equation}
    \label{eq:result-d}
    d(\gamma) = \frac{1}{2t_0^2(at_0 + b)^2}\big((at_0 + b)^4(b - 2at_0)/3 + 2(at_0 + b)^2(\gamma - 2c(\gamma)t_0^3) + (at_0 + b)[t_0f'(t_0) - f(t_0)] - at_0 f(t_0)\big);
\end{equation}
and
\begin{equation}
    \label{eq:result-e}
    e(\gamma)  = \frac{1}{2t_0(at_0 + b)}\left(f(t_0) - (at_0 + b)^4/3 - 2(at_0 + b)(c(\gamma)t_0^3 + d(\gamma)t_0^2 + \gamma)\right). 
\end{equation}
These trisections are all defined over $L$, and each one has a triple singularity at $Q$, as desired.
Here, $a,b$, and $t_0$ are the variables given in the definition of $\famL$.

We will now lay out the verification and derivation of the above family.
This proof is not required to understand the rest of the paper.

A general trisection is of the form:
\begin{equation*}
    T : y = axt + bt + ct^3 + dt^2 + et + h.
\end{equation*}
Note that we are using the variables $a,b$ but they are not \textit{a priori} meant to be the same $a,b$ as in the definition of $\famL$.
If we wish to look at the part of the trisection which is a subset of $\E$, then we are looking at the zero set of the following equation:
\begin{equation*}
    H(x,t) = -(axt + bx + ct^3 + dt^2 + et + h)^2 + x^3 + xf(t) + g(t).
\end{equation*}
Now, let $Q$ be the point $(x_0,y_0,t_0)$; again, $t_0$ is not \textit{a priori} the $t_0$ from the definition of $\famL$.
The trisection $T\subset \E$ has a triple singularity at the point $Q$ if and only if the following six conditions hold:
\begin{align*}
    & H(x_0,t_0) = 0, & \text{(C1)} & &
    & \frac{\del H}{\del t} (x_0,t_0) = 0, & \text{(C2)} & &
    & \frac{\del H}{\del x} (x_0,t_0) = 0, & \text{(C3)} \\
    & \frac{\del^2 H}{\del t^2} (x_0,t_0) = 0, & \text{(C4)} & &
    & \frac{\del^2 H}{\del t \del x}(x_0,t_0) = 0, \text{ and} & \text{(C5)} & &
    & \frac{\del^2 H}{\del x^2}(x_0,t_0) = 0. & \text{(C6)}
\end{align*}
Since $Q\in T$, we already have $y_0 = ax_0t_0 + bx_0 + ct_0^3 + dt_0^2 + et_0 + h$.
Using \magma{}, or another calculator, we see that condition (C6) is equivalent to:
\begin{equation*}
    \text{(C6)} \iff 0 = 2((at_0 + b)^2 - 3x_0)
\end{equation*}
which is satisfied if and only if
\begin{equation}
    \label{eq:x_0}
    x_0 = (at_0 + b)^2/3.
\end{equation}
Plugging in this restriction and solving for (C3) we get:
\begin{equation*}
    \text{(C3)} \iff 0 = (at_0+b)\left[2t_0e + (at_0 + b)^3/3 + 2(ct_0^3 + dt_0^2 + h)\right] - f(t_0), 
\end{equation*}
which can be satisfied by setting
\begin{equation*}
    e = \frac{1}{2t_0(at_0 + b)}\left(f(t_0) - (at_0 + b)^4/3 - 2(at_0 + b)(ct_0^3 + dt_0^2 + h)\right).
\end{equation*}
Substituting this value of $e$ into (C5) yields:
\begin{equation*}
    \text{(C5)} \iff  0 = (at_0 + b)^4(b - 2at_0)/3 + 2(at_0 + b)^2(h - dt_0^2 - 2ct_0^3) + (at_0 + b)[t_0f'(t_0) - f(t_0)] - at_0 f(t_0),
\end{equation*}
which holds by setting
\begin{equation*}
    d = \frac{1}{2t_0^2(at_0 + b)^2}\left((at_0 + b)^4(b - 2at_0)/3 + 2(at_0 + b)^2(h - 2ct_0^3) + (at_0 + b)[t_0f'(t_0) - f(t_0)] - at_0 f(t_0)\right).
\end{equation*}
Further substituting all of the above conditions into (C4) provides:
\begin{multline*}
    \text{(C4)} \iff  0 = - (at_0 + b)^8\left((b - at_0/2)^2 - 7a^2t_0^2/4\right)/9 - (at_0 + b)^7(2h + 2t_0^3 c)/3 \\
    + (at_0 + b)^4\left(c_1(t_0) + t_0^2g''(t_0)\right) - (at_0 + b)^3f(t_0)\left(2h + 2t_0^3 c\right) + c_2(t_0),
\end{multline*}
which is satisfied by setting
\begin{multline*}
    c = \frac{1}{2t_0^3(at_0 + b)^3(f(t_0) + (at_0 + b)^4/3)}\bigg(-(at_0 + b)^8\left((b - at_0/2)^2 - 7a^2t_0^2/4\right)/9 \\
    - 2h(at_0 + b)^7/3 + (at_0 + b)^4\left(c_1(t_0) + t_0^2g''(t_0)\right) - 2h(at_0 + b)^3f(t_0) + c_2(t_0)\bigg),
\end{multline*}
where
\begin{multline}
    \label{eq:c1}
    c_1(t_0) = 5a^2t_0^6f_4 + 3a^2t_0^5f_3 + 5/3a^2t_0^4f_2 + a^2t_0^3f_1 + a^2t_0^2f_0 \\
        + 8abt_0^5f_4 + 13/3abt_0^4f_3 + 2abt_0^3f_2 + abt_0^2f_1 + 4/3abt_0f_0 \\
        + 8/3b^2t_0^4f_4 + b^2t_0^3f_3 - 1/3b^2t_0f_1 
\end{multline}
and
\begin{multline}
    \label{eq:c2}
    c_2(t_0) = a^2t_0^2\bigg( - 13/2t_0^8f_4^2 - 9t_0^7f_3f_4 - 5t_0^6f_2f_4 - 3t_0^6f_3^2 - t_0^5f_1f_4 
    - 3t_0^5f_2f_3 + 3t_0^4f_0f_4 - 1/2t_0^4f_2^2 + 3t_0^3f_0f_3 + t_0^3f_1f_2 \\ 
    + 3t_0^2f_0f_2 + t_0^2f_1^2 + 3t_0f_0f_1 + 3/2f_0^2 \bigg) \\
    + at_0b\bigg(- 17t_0^8f_4^2 - 25t_0^7f_3f_4 - 16t_0^6f_2f_4 - 9t_0^6f_3^2 - 7t_0^5f_1f_4
    - 11t_0^5f_2f_3 + 2t_0^4f_0f_4 - 4t_0^4f_1f_3 - 3t_0^4f_2^2 + 3t_0^3f_0f_3 - t_0^3f_1f_2 \\
    + 4t_0^2f_0f_2 + t_0^2f_1^2 + 5t_0^1f_0f_1 + 3f_0^2\bigg) \\
    + b^2\bigg(- 11t_0^8f_4^2 - 17t_0^7f_3f_4 - 12t_0^6f_2f_4 - 13/2t_0^6f_3^2 - 7t_0^5f_1f_4 
    - 9t_0^5f_2f_3 - 2t_0^4f_0f_4 - 5t_0^4f_1f_3 - 3t_0^4f_2^2 - t_0^3f_0f_3 - 3t_0^3f_1f_2 \\
    - 1/2t_0^2f_1^2 + t_0f_0f_1 + f_0^2\bigg).
\end{multline}
Here the $f_i$s are the $t^i$ coefficients of the polynomial $f(t)$.
Finally, plugging in all of these equations into the remaining conditions gives us:
\begin{equation*}
    \text{(C1)} \iff 0 = \frac{P(t_0)}{(at_0+b)^2}
\end{equation*}
and
\begin{equation*}
    \text{(C2)} \iff 0 = \frac{(at_0+b)P'(t_0) - 2aP(t_0)}{(at_0+b)^3}.
\end{equation*}
We now assume that $a,b,t_0$ are as given by the definition of $\famL$, which ensures that the above two equations equal zero.
The only remaining issue is making sure we are not dividing by zero anywhere.
This is also satisfied by the definition of $\famL$, since each denominator is assumed to be non-zero.

We have a free variable $h$; rewriting it as $\gamma$ and finding $y_0$ from (\ref{eq:x_0}), along with the equation $y_0 = ax_0t_0 + bx_0 + ct_0^3 + dt_0^2 + et_0 + h$, gives us the desired family together with the point $Q$.
$\qed$

\begin{lemma}
    \label{lemma:genus_most_one}
    A trisection with a triple singularity on a del Pezzo surface of degree one has genus at most one.
\end{lemma}

\proof

Use (\ref{eq:genus_formula}); plugging in the triple singularity gives a reduction of 3, for a maximum genus of one.
$\qed$

\begin{proposition}
    \label{prop:unique}
    Suppose $S\in\famL$, and suppose that $\T_\gamma$ is the family obtained in Lemma \ref{lemma:trisection}.
    Let $R = (x_R,y_R,t_R)\in\E$ with $t_R\neq t_0$, the $t$-coordinate of $Q$.
    Then, there exists a unique choice of $\gamma\in L(x_R,y_R,t_R)$ such that $R\in \T_\gamma$.
    In particular, if $R\in\E$ is $L$-rational, then there is a $\gamma\in L$ such that $R\in \T_\gamma/L$.
\end{proposition}

\proof

The family of trisections, $\T_\gamma$, is defined by equations (\ref{eq:result-c}), (\ref{eq:result-d}), and (\ref{eq:result-e}) from Lemma \ref{lemma:trisection}.
Plugging in all the restrictions, we can use \magma{} to calculate that the equation of definition of $\T_\gamma$ can be written in the form:
\begin{equation*}
    \T_\gamma : y = (1-t/t_0)^3 \gamma + A(x,t),
\end{equation*}
where $A(x,t)$ is polynomial in $x,t$ with coefficients in $L$, all of which are polynomial in $a,b,t_0$ and the coefficients of $f,g$.
So, $R\in\T_\gamma$ if and only if:
\begin{equation*}
    y_R = (1 - t_R/t_0)^3 \gamma + A(x_R,t_R).
\end{equation*}
Furthermore, since $t_R\neq t_0$, we can isolate for $\gamma$, hence
\begin{equation*}
    \gamma = \frac{y_R - A(x_R,t_R)}{(1-t_R/t_0)^3}.
\end{equation*}
From this we see both that $\gamma$ exists and is unique.
Since $A$ is a polynomial in $x,t$ over $L$, and $t_0\in L$ we see that $\gamma\in L(x_R,y_R,t_R)$.
$\qed$

\begin{definition}
    \label{def:T-R}
    Supposing we have the same setup as Proposition \ref{prop:unique}, we can define $\T_R$ to be the unique member of the family $\T_\gamma$ with $R\in\T_R$.
\end{definition}

\begin{lemma}
    \label{lemma:three-choices}
    Suppose $R=(x_R,y_R,t_R)\in\E$, with $t_R\neq t_0$, and let $P=(x_P,y_P,t_P)\in\E$ be a point different from $R$ but with the same $t$ value.
    Then, there are at most two such points $P\in \T_R$.
\end{lemma}

\proof 

Let $t_R = t_P = t_1$; then, both $R,P\in\E_{t_1}$.
Suppose that $P\in \T_R$; this means that $\T_R = \T_P$, so, in particular, they both belong to $\T_\gamma$ for the same choice of $\gamma$.
Through calculation in \magma{}, or just by looking at equation (\ref{eq:trisection}), we can see that $A(x,t)$ is linear in $x$, with coefficients which are polynomials in $t$.
Thus, $y_R$ can be derived from $x_R,t_R$ from equation (\ref{eq:trisection}), and the same is true for $y_P$.
Plugging these equations back into equation (\ref{eq:affinesurface}), the equation of definition of $\E$, we get equations which are cubic in $x_R$ and $x_P$ with coefficients which are polynomial in $t_R$ and $t_P$, respectively.
However, $t_R = t_P = t_1$, so these are, in fact, the same equations, but with $x_R$ replaced with $x_P$; hence, $x_P$ and $x_R$ both satisfy the same cubic polynomial in $x$.
Therefore there are at most two solutions to $x_P$ different from $x_R$, which gives two possible choices of $P$.
$\qed$

\begin{theorem}[Theorem \ref{thm:intro-main}]
    \label{thm:main}
    Let $L/k$ be a field extension and $S\in\famL$ a surface.
    Further assume that $S$ contains a fibre $\F : Z = t_1 W$, with $Q\nin \F$, such that $|\F(L)| = \infty$.
    Then the family from Lemma \ref{lemma:trisection} is a pencil in $S$ and each element of the pencil is a curve of genus zero or genus one defined over $L$.
\end{theorem}

\proof

Combine Lemma \ref{lemma:trisection} and Lemma \ref{lemma:genus_most_one} to show the genus of each element of the family $\T_\gamma$.
To see that the family $\T_\gamma$ is a pencil, we must show that it is, in fact, one-dimensional.
As seen in \ref{lemma:trisection}, this family has a single free variable, so it is enough to show that it contains infinitely many distinct curves defined over $L$.
This follows from Proposition \ref{prop:unique} and Lemma \ref{lemma:three-choices}, since each $R\in|\F(L)|$ defines a curve $\T_R$ in the family, and the mapping $R\mapsto \T_R$ is, at worst, three-to-one.
Since $|\F(L)| = \infty$, we see that there are infinitely many distinct $\T_\gamma$ defined over $L$.
$\qed$
\section{A Rational Point on $\T_R$}
\label{sec:rational-point}

\begin{proposition}
    \label{prop:rational-point}
    Suppose that $R\in\E$ with $R = (x_R,y_R,t_R)$ and $t_R\neq t_0$, and further suppose that $\T_R$ has genus one.
    Let $P_1,P_2$ be the two other points on the $t_R$ fibre of $\T_R$.
    Then, $P_1 + P_2$ is $L(x_R,y_R,t_R)$-rational and, furthermore, $R,P_1,P_2$ are all distinct.
\end{proposition}

\proof

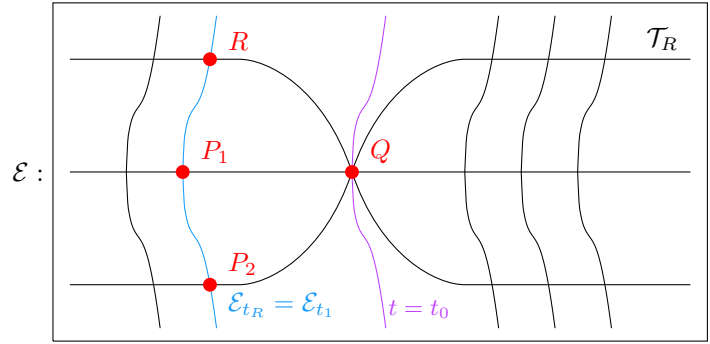
\begin{figure}[H]
    \centering
    \begin{tikzpicture}[scale=1.5]
        \draw (-1.5,0) -- (0,0) .. controls (0.25,0) and (0.75,0.25) .. (1,1) .. controls (1.25,1.75) and (1.75,2) .. (2,2) -- (4,2);
        \draw (-1.5,1) -- (0,1) -- (1,1) -- (2,1) -- (4,1);
        \draw (-1.5,2) -- (0,2) .. controls (0.25,2) and (0.75,1.75) .. (1,1) .. controls (1.25,0.25) and (1.75,0) .. (2,0) -- (4,0);

        \draw[variable=\t, domain=-1:2,smooth,samples=\sample,color=functioncolor] plot ({(\xscale)*(\t+1)+1},{(\yscale)*(\t^3 + \t + 2)^(0.5)+1});
        \draw[variable=\t, domain=-1:2,smooth,samples=\sample,color=functioncolor] plot ({(\xscale)*(\t+1)+1},{-(\yscale)*(\t^3 + \t + 2)^(0.5)+1});

        \draw[variable=\t, domain=-1:2,smooth,samples=\sample] plot ({(\xscale)*(\t+1)+3},{(\yscale)*(\t^3 + \t + 2)^(0.5)+1});
        \draw[variable=\t, domain=-1:2,smooth,samples=\sample] plot ({(\xscale)*(\t+1)+3},{-(\yscale)*(\t^3 + \t + 2)^(0.5)+1});

        \draw[variable=\t, domain=-1:2,smooth,samples=\sample] plot ({(\xscale)*(\t+1)+2.5},{(\yscale)*(\t^3 + \t + 2)^(0.5)+1});
        \draw[variable=\t, domain=-1:2,smooth,samples=\sample] plot ({(\xscale)*(\t+1)+2.5},{-(\yscale)*(\t^3 + \t + 2)^(0.5)+1});

        \draw[variable=\t, domain=-1:2,smooth,samples=\sample] plot ({(\xscale)*(\t+1)+2},{(\yscale)*(\t^3 + \t + 2)^(0.5)+1});
        \draw[variable=\t, domain=-1:2,smooth,samples=\sample] plot ({(\xscale)*(\t+1)+2},{-(\yscale)*(\t^3 + \t + 2)^(0.5)+1});

        \draw[variable=\t, domain=-1:2,smooth,samples=\sample] plot ({(\xscale)*(\t+1)-1},{(\yscale)*(\t^3 + \t + 2)^(0.5)+1});
        \draw[variable=\t, domain=-1:2,smooth,samples=\sample] plot ({(\xscale)*(\t+1)-1},{-(\yscale)*(\t^3 + \t + 2)^(0.5)+1});

        \draw[variable=\t, domain=-1:2,smooth,samples=\sample,color=fibrecolor] plot ({(\xscale)*(\t+1)-0.5},{(\yscale)*(\t^3 + \t + 2)^(0.5)+1});
        \draw[variable=\t, domain=-1:2,smooth,samples=\sample,color=fibrecolor] plot ({(\xscale)*(\t+1)-0.5},{-(\yscale)*(\t^3 + \t + 2)^(0.5)+1});

        \fill[highlighted] (1,1) circle (1.75 pt);
        \fill[highlighted] (-0.5,1) circle (1.75 pt);
        \fill[highlighted] (-0.5+\xoff,0) circle (1.75 pt);
        \fill[highlighted] (-0.5+\xoff,2) circle (1.75 pt);

        \draw (-1.65,-0.5) -- (4.15,-0.5) -- (4.15,2.5) -- (-1.65,2.5) -- (-1.65,-0.5);

        \draw (-1.65,1) node [anchor=east] {$\E:$};
        \draw (3.75,2) node [anchor=south] {$\T_R$};
        \draw[color=highlighted] (1,1) node [anchor=south west] {\hspace{0 pt} $Q$};
        \draw[color=highlighted] (-0.5,1) node [anchor=south west] {\hspace{0 pt} $P_1$};
        \draw[color=highlighted] (-0.5+\xoff,2) node [anchor=south west] {\hspace{0 pt} $R$};
        \draw[color=highlighted] (-0.5+\xoff,0) node [anchor=south west] {\hspace{0 pt} $P_2$};
        \draw[color=fibrecolor] (-0.5+\xoff,-0.375) node [anchor=south west] {\hspace{0 pt} $\E_{t_R} = \E_{t_1}$};
        \draw[color=functioncolor] (1+\xoff,-0.375) node [anchor=south west] {\hspace{0 pt}{\small $t = t_0$}};
    \end{tikzpicture}
    \caption{An example of a trisection $\T_R$ in the family $\{\T_R\}_{R\in\E_{t_1}}$. 
    Fibres containing the relevant points have been highlighted.}
    \label{fig:whatsgoingon}
\end{figure}

Let $R\in\E$ be an arbitrary point with $R = (x_R,y_R,t_R)$ and $t_R \neq t_0$.
Now, consider the trisection $\T_R$ and the fibre $\E_{t_R}$; by the proof of Lemma \ref{lemma:three-choices}, we know that this trisection intersects the fibre at three points, counting multiplicity.
One of these points is $R$; call the other two $P_1,P_2$.
Figure \ref{fig:whatsgoingon} shows the general picture of the trisection.

We know $R,P_1,P_2$ are all distinct, since if any two were the same point, then the trisection would have another singularity, and so by equation (\ref{eq:genus_formula}), the genus of $\T_R$ would be zero instead of one.
The points $P_1,P_2$ are both on the fibre $\E_{t_R}$, so they both have $t$-coordinate $t_R$.
Similarly, since they are both on the trisection $\T_R$, both points' $y$-coordinates can be deduced from their $t$- and $x$-coordinates.
We see that the two points are the solutions to:
\begin{equation*}
    \left(axt_R + bx + ct_R^3 + dt_R^2 + et_R + h\right)^2 = x^3 + f(t_R)x + g(t_R), 
\end{equation*}
which is a cubic polynomial in $x$.
This has a root at $x_R$ since $R\in\T_R$; the other two roots are the $x$-coordinates of $P_1$ and $P_2$.
Suppose that $L'/L(x_R,y_R,t_R)$ is the field extension of $L(x_R,y_R,t_R)$ which contains these two additional roots; then, the points $P_1$ and $P_2$ are defined over $L'$ since the $t$-coordinate is just $t_R$, and the $y$-coordinate is polynomial over $L'$ in its $x$- and $t$-coordinates.

Now, look at the Galois action of $L'/L(x_R,y_R,t_R)$ acting on the points of $P_1,P_2$ by acting on their coordinates.
Clearly, it will map the set $\{P_1,P_2\}$ into itself, either because $P_1$, and hence also $P_2$, is $L(x_R,y_R,t_R)$-rational, or because the two $x$-coordinates are two non-$L(x_R,y_R,t_R)$-rational roots of the polynomial of which $L'/L(x_R,y_R,t_R)$ is the splitting field.
In the second case, the two non-$L(x_R,y_R,t_R)$-rational roots will be mapped to one another.

Looking at the curve $\T_R$, which has geometric genus one, by construction, it is a singular curve.
It has a triple singularity at $Q$ and is non-singular everywhere else, by the same argument that shows $R,P_1,P_2$ are all distinct.
Let $E$ be the projective closure of the desingularization of $\T_R$ at $Q$; then, $E$ is a non-singular projective curve of geometric genus one.
This desingularization does not change the field of definition, so $E$ is defined over $L(x_R,y_R,t_R)$.
The point $R$ has a corresponding point on $E$, which we will call $\hat R$; it still has the same coordinates.
This means since $R$ is an $L(x_R,y_R,t_R)$-rational point, so too is $\hat R$.
Now we look at the points $\hat P_1,\hat P_2$, and the Galois action of $L'/L(x_R,y_R,t_R)$; it will map the set $\{\hat P_1,\hat P_2\}$ to itself.
The Galois action also commutes with the group action on the elliptic curve $(E,\hat R)$, since the group action can be expressed using rational functions over $L(x_R,y_R,t_R)$.
So, $\hat P_1 + \hat P_2$ must therefore be mapped to itself by the Galois action, meaning that $\hat P_1 + \hat P_2$ is an $L(x_R,y_R,t_R)$-rational point on $(E,\hat R)$, with a corresponding $L(x_R,y_R,t_R)$-rational point on $\T_R$.
$\qed$
\section{Next Steps and Applications}
\label{sec:order}
\subsection{Zariski Density}
\label{sec:Zariski}

In Section \ref{sec:rational-point} we found an $L(x_R,y_R,z_R)$-rational point $P_R := \hat P_1 + \hat P_2$ on $(E,\hat R)$.
If $L'/L$ is an arbitrary field extension, and $P$ is an $L'$-rational point on a elliptic curve $\E/L$, then $nP$ is also an $L'$-rational point in $\E$ for any $n\in\N$.
This holds because the formula for adding points on a rational curve is a rational function of the coefficients of the points being added.
Hence we know that $nP_R$ is also $L(x_R,y_R,z_R)$-rational for any $n\in\N$.
In particular, we wish to show that the order of $P_R$ is always infinite.
If this were the case, we could use the following theorem to prove Zariski density for any surface $S\in\famL$; however, to show this we still need to find a Weierstrass model for each trisection in the family of trisections $\T_R$.

\begin{theorem}[Theorem \ref{thm:intro-ZD}]
    \label{thm:ZD}
    Let $L/k$ be a field extension, and suppose $S\in\famL$, with a fibre $\F\subset S$, such that $|\F(L)| = \infty$.
    Further, suppose that one of the following three hold:
    \begin{enumerate}[label=(\roman*)]
        \item there exists an $R_0\in \F(L)$ such that $\T_{R_0}$ has genus zero, or;

        \item for infinitely many $R\in\F(L)$, $\T_R$ has genus one and the order of $P_R$ is infinite in $\T_R$, or;

        \item there exists an $R^*\in\F(L)$, with $\T_{R^*}$ having genus one, $P_{R^*}$ having infinite order, and for infinitely many $n\in\N$ the point $nP_{R^*}$ has infinite order on its fibre.
    \end{enumerate}
    Then, the $L$-rational points, $S(L)$, are Zariski dense in $S$. 
\end{theorem}

\proof

Assume hypothesis (i).
Note that by (\ref{eq:genus_formula}) and the construction of the pencil $\{T_R\}_{R\in\F}$, the curve $\T_{R_0}$ has one singularity different from $Q$, which cannot be a triple singularity, since the genus is zero and cannot be negative.
We look at two cases.
Either the additional singularity occurs at the point $R_0$, in which case $R_0$ is a double point and by Lemma \ref{lemma:three-choices}, there is a unique non-singular point $P\in\T_{R_0}\cap\F$ which will be an $L$-rational point since $R_0$ is a double point.
Otherwise, $R_0$ is a non-singular $L$-rational point.
In either case, we see that the genus zero curve $\T_{R_0}$ contains a non-singular $L$-rational point.
Hence, by the same argument as used in \cite[Theorem 1.1]{DW_2022}, we see that the $L$-rational points of $S$ are Zariski dense.

Now moving onto hypothesis (ii); suppose that the order of $P_R$ is infinite, for infinitely many choices of $R\in\F$.
In this case, let $R\in\F(L)$, and look at $\T_R$, which we assume has genus one.
Then, for each $\T_R$, by Proposition \ref{prop:rational-point}, there is an $L$-rational point $P_R\in\T_R$, which, by our assumption, for infinitely many of the $R$'s, will have infinite order.
Furthermore, except for at most two choices of $n$, the points $nP_R\nin\F$ by Lemma \ref{lemma:three-choices}.
By Proposition \ref{prop:unique}, each point $nP_R$ can only occur on a single member $\T_R$; by construction, it will belong to the same $\T_R$ as $P_R$.
Therefore there are infinitely many distinct genus one curves $\T_R$, which all contain a point $P_R$, which has infinite order on the corresponding elliptic curve.
I.e., there are infinitely many distinct curves defined over $L$, each of which has infinitely many distinct $L$-rational points.
Therefore, the $L$-rational points of $S$ cannot be contained in any proper closed subset of $S$, and so the $L$-rational points must be dense in the Zariski topology.

Finally we have hypothesis (iii); in this case, we assume there is a single point $R^*$ such that there are infinitely many points $nP_{R^*}\in\T_{R^*}$ which are non-torsion on their fibre.
Then, we see that similar to the case of hypothesis (ii), we have infinitely many distinct curves which all contain infinitely many distinct $L$-rational points.
In particular, these curves are the fibres of each $nP_{R^*}$ for which $nP_{R^*}$ is non-torsion.
Hence, the $L$-rational points cannot be contained in a single proper closed subset of $S$ and so must be dense in the Zariski topology.
$\qed$

The standard way to prove that a point has infinite order requires a Weierstrass model for the elliptic curve, which is possible to find using \magma{}; however, in the general setting of Theorem \ref{thm:ZD}, we ran into issues of extensive computation time.
We were able to find the model for several subfamilies of surfaces and these families are shown in Proposition \ref{prop:models}.

Case (i) of Theorem \ref{thm:intro-ZD} requires that the surface has a trisection of genus zero.
Theorem \ref{thm:famone} provides a large family of surfaces with this property.

\begin{definition}
    \label{def:famone}
    Fix a field $k$ of characteristic zero.
    Let $\famone$ be the family consisting of all rational elliptic surfaces of the form:
    \begin{equation*}
        \E : y^2 = x^3 + f(t) x + g(t),
    \end{equation*}
    with $f(t) = f_0 + f_1t + f_2t^2 + f_3t^3 + f_4t^4$, and $g(t) = G_0 + G_1t + G_2 t^2 + G_3t^3 + g_4t^4 + g_5t^5 + g_6t^6$.
    Here $f_i,g_i\in k$ are all free variables, and the variables $G_i$ are all rational functions of the coefficients $f_i,g_i$ along with five other free variables, $a,b,t_0,t_1,x_1\in k$.
    These additional free variables, along with seven other rational functions $X_0,Y_0,Y_1,C,D,E,H$ of the free variables $f_i,g_i,a,b,t_0,t_1,x_1$, are used in Theorem \ref{thm:famone} to define a trisection on the surface $\E$.
    In particular, the equations for $G_i$ and the other rational functions are found \thecode.
\end{definition}

\begin{theorem}[Theorem \ref{thm:intro-DP1}]
    \label{thm:famone}
    Let $\E\in\famone$.
    Then, $\E$ contains the trisection $T$ defined by:
    \begin{equation*}
        T : axt + bt + Ct^3 + Dt^2 + Et + H,
    \end{equation*}
    which has a triple singularity at the point $Q = (X_0,Y_0,t_0)$ and a double singularity at the point $R = (x_1, Y_1, t_1)$.
    Recall the definitions of the rational functions $X_0,Y_0,Y_1,C,D,E,H$ from Definition \ref{def:famone}.
    On $\E$, the $k$-rational points are Zariski dense.
\end{theorem}

\begin{remark}
    As stated in Definition \ref{def:famone}, the variables $a,b,t_0,t_1,x_1,f_i,g_i\in k$ are all free variables.
    In particular, $\famone$ consists of thirteen free variables.
    These rational elliptic surfaces $\E$ generically have only irreducible bad fibres, and hence are the blow-up of a del Pezzo surface of degree one.
    The corresponding del Pezzo surface $S$ is a member of $\famLk$.
\end{remark}

\proof

We start with a generic elliptic surface of the form
\begin{equation*}
    \E : y^2 = x^3 + f(t)x + g(t),
\end{equation*}
with $f(t),g(t)\in k[t]$ polynomials of degree four and six, respectively.
Let them have coefficients $f_i$ and $g_j$ with $0 \leq i \leq 4$ and $0 \leq j \leq 6$.
Let $T$ be a trisection given by:
\begin{equation*}
    T : y = axt + bx + ct^3 + dt^2 + et + h,
\end{equation*}
where $a,b,c,d,e,h\in k$ are free variables.
Let $Q = (x_0,y_0,t_0)$, where $x_0,y_0,t_0\in k$ are free variables, and suppose that $t_0(at_0 + b) \neq 0$ and $f(t_0) \neq -(at_0 + b)^4/3$.
Setting $x_0,c,d,e$ to (\ref{eq:x_0}), (\ref{eq:result-c}), (\ref{eq:result-d}), and (\ref{eq:result-e}), respectively, with free variable $\gamma = h$, and setting $y_0$ as in Lemma \ref{lemma:trisection}, we obtain a trisection $T$ defined over $k$, and a $k$-rational point $Q = ((at_0 + b)^2/3,(at_0 + b)^3/6 + f(t_0)/(2(at_0 + b)),t_0)$.

This point is a triple singularity on the trisection $T$, so long as the polynomial $P(t)$ from equation (\ref{eq:condition}) has a double root at $t = t_0$.
In this case the corresponding del Pezzo surface $S$ will be in $\famLk$.
To force this condition on $P(t)$ to hold, we can set $P(t_0) = 0$ and solve for $g_0$, from which we obtain:
\begin{equation*}
    g(t_0) = \frac{1}{(at_0 + b)^2}\left(\frac{f(t_0)^2}{4} - \frac{(at_0 + b)^4f(t_0)}{6} - \frac{(at_0 + b)^8}{108}\right),
\end{equation*}
along with:
\begin{equation*}
    g_0 = \frac{1}{(at_0 + b)^2}\bigg(\frac{f(t_0)^2}{4} - \frac{(at_0 + b)^4f(t_0)}{6} - \frac{(at_0 + b)^8}{108} - (at_0 + b)^2\left(g_6t_0^6 + g_5t_0^5 + g_4t_0^4 + g_3t_0^3 + g_2t_0^2 + g_1t_0\right)\bigg).
\end{equation*}
We also need $P'(t_0) = 0$; since we have $g(t_0)$ from above, we have the following equation for $(at_0 + b)P'(t_0)$:
\begin{align*}
    (at_0 + b)P'(t_0) & =\frac{-(at_0 + b)f(t_0)f'(t_0)}{2} + \frac{2a(at_0 + b)^4f(t_0)}{3} + \frac{f'(t_0)(at_0 + b)^5}{6} + 2a(at_0 + b)^2g(t_0) \\
    & \hspace{10 em} + (at_0 + b)^3g'(t_0) + \frac{2a(at_0 + b)^8}{27} & \\
    & =\frac{af(t_0)^2}{2} - \frac{(at_0 + b)f(t_0)f'(t_0)}{2} + \frac{a(at_0 + b)^4f(t_0)}{3} + \frac{f'(t_0)(at_0 + b)^5}{6} & \\
    & \hspace{10 em} + (at_0 + b)^3g'(t_0) + \frac{a(at_0 + b)^8}{18}. &
\end{align*}
Hence, we get the following equation for $g_1$:
\begin{multline*}
    g_1 = \frac{1}{(at_0 + b)^3}\bigg(-\frac{af(t_0)^2}{2} + \frac{(at_0 + b)f(t_0)f'(t_0)}{2} - \frac{a(at_0 + b)^4f(t_0)}{3} - \frac{f'(t_0)(at_0 + b)^5}{6} - \frac{a(at_0 + b)^8}{18} \\
    - (at_0 + b)^3\left(6g_6t_0^5 + 5g_5t_0^4 + 4g_4t_0^3 + 3g_3t_0^2 + 2g_2t_0\right)\bigg).
\end{multline*}
From here, we add in a new point $R = (x_1,y_1,t_1)\in\E$; once again, we have conditions for the trisection to go through $R$ with multiplicity two.
They are the following (where $H(x,t)$ is the same as from Lemma \ref{lemma:trisection}):
\begin{align*}
    & H(x_0,t_0) = 0, & \text{(C7)} & &
    & \frac{\del H}{\del t} (x_0,t_0) = 0, \text{ and} & \text{(C8)} & &
    & \frac{\del H}{\del x} (x_0,t_0) = 0. & \text{(C9)} 
\end{align*}
By using \magma{}, we see that (C9) is satisfied when
\begin{multline*}
    h = \frac{1}{2(at_0 + b)^3(at_1 + b)(t_0 - t_1)^3((at_0 + b)^4 + 3f(t_0))}\bigg(-6t_1t_0^2(t_0 - t_1)^2(at_0 + b)^4(at_1 + b)\bigg[g_2  + 3g_3t_0  + 6g_4t_0^2 \\
    + 10g_5t_0^3 + 15g_6t_0^4\bigg] + (at_0 + b)^7\bigg[x_1t_0^3(3x_1 - 2(at_1 + b)^2) + t_1(at_0 + b)(at_1 + b)\bigg(-5/6a^2t_0^4 + 5/3a^2t_0^3t_1 \\
    - 1/2a^2t_0^2t_1^2 + abt_0^2t_1 - 1/3abt_0t_1^2 + b^2t_0^2 - b^2t_0t_1 + 1/3b^2t_1^2\bigg)\bigg] - 6x_1t_0^3(at_0 + b)^3\bigg[(at_1 + b)^2 - 3/2x_1\bigg]\bigg[f(t_0) \\
    + t_0(at_0 + b)^4(t_0^2(at_0 + b)^3f_0 + t_1h_1) + (t_0 - t_1)^2h_2\bigg]\bigg),
\end{multline*}
where
\begin{multline*}
    h_1 = a^3t_0\bigg(-15t_0^6t_1f_4 + 30t_0^5t_1^2f_4 - 10t_0^5t_1f_3 - 14t_0^4t_1^3f_4 + 20t_0^4t_1^2f_3 - 6t_0^4t_1f_2 + t_0^4f_1 - 9t_0^3t_1^3f_3 + 12t_0^3t_1^2f_2 - 6t_0^3t_1f_1 \\
    - 5t_0^2t_1^3f_2 + 9t_0^2t_1^2f_1 - 7t_0^2t_1f_0 - 3t_0t_1^3f_1 + 10t_0t_1^2f_0 - 3t_1^3f_0\bigg) \\ 
    + a^2b\bigg(-15t_0^7f_4 + 10t_0^6t_1f_4 - 10t_0^6f_3 + 29t_0^5t_1^2f_4 + 8t_0^5t_1f_3 - 7t_0^5f_2 - 21t_0^4t_1^3f_4 + 18t_0^4t_1^2f_3 + 9t_0^4t_1f_2 - 3t_0^4f_1 - 13t_0^3t_1^3f_3 + 7t_0^3t_1^2f_2 \\
    + 4t_0^3t_1f_1 - 7t_0^3f_0 - 6t_0^2t_1^3f_2 + 5t_0^2t_1^2f_1 + 2t_0^2t_1f_0 - 3t_0t_1^3f_1 + 9t_0t_1^2f_0 - 4t_1^3f_0\bigg) \\
    + ab^2\bigg(-20t_0^6f_4 + 40t_0^5t_1f_4 - 11t_0^5f_3 - 12t_0^4t_1^2f_4 + 24t_0^4t_1f_3 - 6t_0^4f_2 - 5t_0^3t_1^3f_4 - 7t_0^3t_1^2f_3 + 17t_0^3t_1f_2 - 2t_0^3f_1 - 3t_0^2t_1^3f_3 - 8t_0^2t_1^2f_2 \\
    + 10t_0^2t_1f_1 - 8t_0^2f_0 - 6t_0t_1^2f_1 + 12t_0t_1f_0 + t_1^3f_1 - 4t_1^2f_0\bigg) \\
    + b^3\bigg(-4t_0^5f_4 + 12t_0^4t_1f_4 - 8t_0^3t_1^2f_4 + 3t_0^3t_1f_3 + 2t_0^3f_2 +
    t_0^2t_1^3f_4 - 2t_0^2t_1^2f_3 - t_0^2t_1f_2 + 3t_0^2f_1 - 3t_0t_1f_1 + t_1^2f_1\bigg),
\end{multline*}
and
\begin{multline*}
    h_2 = a^3t_0^2\bigg(45/2t_0^8t_1^2f_4^2 + 30t_0^7t_1^2f_3f_4 + 18t_0^6t_1^2f_2f_4 + 9t_0^6t_1^2f_3^2 + 3t_0^6t_1f_1f_4 + 3t_0^6f_0f_4 + 6t_0^5t_1^2f_1f_4 + 9t_0^5t_1^2f_2f_3 + 6t_0^5t_1f_0f_4 \\
    + 3t_0^5t_1f_1f_3 + 3t_0^5f_0f_3 - 6t_0^4t_1^2f_0f_4 + 3/2t_0^4t_1^2f_2^2 + 6t_0^4t_1f_0f_3 + 3t_0^4t_1f_1f_2 + 3t_0^4f_0f_2 - 9t_0^3t_1^2f_0f_3 - 3t_0^3t_1^2f_1f_2 \\ 
    + 6t_0^3t_1f_0f_2 + 3t_0^3t_1f_1^2 + 3t_0^3f_0f_1 - 9t_0^2t_1^2f_0f_2 - 3t_0^2t_1^2f_1^2 + 9t_0^2t_1f_0f_1 + 3t_0^2f_0^2 - 9t_0t_1^2f_0f_1 + 6t_0t_1f_0^2 - 9/2t_1^2f_0^2\bigg) \\
    + a^2bt_0\bigg(45/2t_0^9t_1f_4^2 + 60t_0^8t_1^2f_4^2 + 30t_0^8t_1f_3f_4 + 84t_0^7t_1^2f_3f_4 + 15t_0^7t_1f_2f_4 + 9t_0^7t_1f_3^2 + 57t_0^6t_1^2f_2f_4 + 27t_0^6t_1^2f_3^2 + 9t_0^6t_1f_1f_4 \\
    + 6t_0^6t_1f_2f_3 + 9t_0^6f_0f_4 + 30t_0^5t_1^2f_1f_4 + 33t_0^5t_1^2f_2f_3 + 3t_0^5t_1f_0f_4 + 3t_0^5t_1f_1f_3 - 3/2t_0^5t_1f_2^2 + 9t_0^5f_0f_3 + 3t_0^4t_1^2f_0f_4 \\
    + 12t_0^4t_1^2f_1f_3 + 9t_0^4t_1^2f_2^2 - 3t_0^4t_1f_1f_2 + 9t_0^4f_0f_2 - 9t_0^3t_1^2f_0f_3 + 3t_0^3t_1^2f_1f_2 - 3t_0^3t_1f_0f_2 + 9t_0^3f_0f_1 - 12t_0^2t_1^2f_0f_2 \\ 
    - 3t_0^2t_1^2f_1^2 + 3t_0^2t_1f_0f_1 + 9t_0^2f_0^2 - 15t_0t_1^2f_0f_1 + 9/2t_0t_1f_0^2 - 9t_1^2f_0^2\bigg) \\
    + ab^2\bigg(60t_0^9t_1f_4^2 + 42t_0^8t_1^2f_4^2 + 87t_0^8t_1f_3f_4 + 60t_0^7t_1^2f_3f_4 + 54t_0^7t_1f_2f_4 + 30t_0^7t_1f_3^2 + 45t_0^6t_1^2f_2f_4 + 39/2t_0^6t_1^2f_3^2 + 30t_0^6t_1f_1f_4 \\
    + 33t_0^6t_1f_2f_3 + 9t_0^6f_0f_4 + 30t_0^5t_1^2f_1f_4 + 27t_0^5t_1^2f_2f_3 + 6t_0^5t_1f_0f_4 + 15t_0^5t_1f_1f_3 + 6t_0^5t_1f_2^2 + 9t_0^5f_0f_3 + 15t_0^4t_1^2f_0f_4 \\
    + 15t_0^4t_1^2f_1f_3 + 9t_0^4t_1^2f_2^2 - 3t_0^4t_1f_0f_3 + 9t_0^4f_0f_2 + 3t_0^3t_1^2f_0f_3 + 9t_0^3t_1^2f_1f_2 - 12t_0^3t_1f_0f_2 - 3t_0^3t_1f_1^2 + 9t_0^3f_0f_1 \\ 
    + 3/2t_0^2t_1^2f_1^2 - 12t_0^2t_1f_0f_1 + 9t_0^2f_0^2 - 3t_0t_1^2f_0f_1 - 6t_0t_1f_0^2 - 3t_1^2f_0^2\bigg) \\
    + b^3\bigg(39t_0^8t_1f_4^2 + 3t_0^7t_1^2f_4^2 + 60t_0^7t_1f_3f_4 + 3t_0^6t_1^2f_3f_4 + 42t_0^6t_1f_2f_4 + 45/2t_0^6t_1f_3^2 + 3t_0^5t_1^2f_2f_4 + 27t_0^5t_1f_1f_4 + 30t_0^5t_1f_2f_3 \\ 
    + 3t_0^5f_0f_4 + 3t_0^4t_1^2f_1f_4 + 12t_0^4t_1f_0f_4 + 18t_0^4t_1f_1f_3 + 9t_0^4t_1f_2^2 + 3t_0^4f_0f_3 + 3t_0^3t_1^2f_0f_4 + 6t_0^3t_1f_0f_3 + 9t_0^3t_1f_1f_2 \\
    + 3t_0^3f_0f_2 + 3/2t_0^2t_1f_1^2 + 3t_0^2f_0f_1 - 3t_0t_1f_0f_1 + 3t_0f_0^2 - 3t_1f_0^2\bigg).
\end{multline*}
We also find that (C7) is satisfied when 
\begin{multline*}
    g_2 = -\frac{1}{(at_0 + b)^3(at_1 + b)^2(t_0 - t_1)^2}\bigg((at_0 + b)^4(at_1 + b)^2\bigg[a\bigg(f_4(5/6t_0^5 - t_0^4t_1) + f_3(2/3t_0^4 - 5/6t_0^3t_1) +
    f_2(1/2t_0^3 - 2/3t_0^2t_1) \\ + f_1(1/3t_0^2 - 1/2t_0t_1) + f_0(1/6t_0 - 1/3t_1)\bigg) + 
    b\bigg(f_4(1/2t_0^4 - 2/3t_0^3t_1) + f_3(1/3t_0^3 - 1/2t_0^2t_1) + f_2(1/6t_0^2 - 1/3t_0t_1) \\
    - 1/6f_1t_1 - 1/6f_0\bigg)\bigg] + 
    (t_0 - t_1)^2G + x_1(at_0 + b)^3\bigg[(at_1 + b)^2(x_1^2 + f(t_1)) - (9/4x_1^3 + 3/2x_1f(t_1))\bigg] \\
    +
    (t_0 - t_1)^2(at_0 + b)^3(at_1 + b)^2\bigg[g_6(5t_0^4 + 4t_0^3t_1 + 3t_0^2t_1^2 + 2t_0t_1^3 + t_1^4) + g_5(4t_0^3 + 3t_0^2t_1 +
    2t_0t_1^2 + t_1^3) + g_4(3t_0^2 + 2t_0t_1 + t_1^2) \\
    + g_3(2t_0 + t_1)\bigg] - (at_0 + b)^8(at_1 + b)^2\bigg[1/18at_1 - 5/108at_0 + 1/108b\bigg]\bigg),
\end{multline*}
where
\begin{multline*}
    G = a^3\bigg(-5/4t_0^7t_1^2f_4^2 - t_0^6t_1^3f_4^2 - 2t_0^6t_1^2f_4f_3 - 3/4t_0^5t_1^4f_4^2 -
    3/2t_0^5t_1^3f_4f_3 - 3/2t_0^5t_1^2f_4f_2 - 3/4t_0^5t_1^2f_3^2 -
    1/2t_0^4t_1^5f_4^2 - t_0^4t_1^4f_4f_3 \\
    - t_0^4t_1^3f_4f_2 - 1/2t_0^4t_1^3f_3^2
    - t_0^4t_1^2f_4f_1 - t_0^4t_1^2f_3f_2 - 1/4t_0^3t_1^6f_4^2 -
    1/2t_0^3t_1^5f_4f_3 - 1/2t_0^3t_1^4f_4f_2 - 1/4t_0^3t_1^4f_3^2 -
    1/2t_0^3t_1^3f_4f_1 \\
    - 1/2t_0^3t_1^3f_3f_2 - 1/2t_0^3t_1^2f_4f_0 -
    1/2t_0^3t_1^2f_3f_1 - 1/4t_0^3t_1^2f_2^2 - 1/2t_0t_1f_1f_0 - 1/4t_0f_0^2 -
    1/2t_1f_0^2\bigg) \\
    + a^2b\bigg(-5/2t_0^7t_1f_4^2 - 15/4t_0^6t_1^2f_4^2 - 4t_0^6t_1f_4f_3 - 3t_0^5t_1^3f_4^2 -
    6t_0^5t_1^2f_4f_3 - 3t_0^5t_1f_4f_2 - 3/2t_0^5t_1f_3^2 - 9/4t_0^4t_1^4f_4^2
    - 9/2t_0^4t_1^3f_4f_3 \\
    - 9/2t_0^4t_1^2f_4f_2 - 9/4t_0^4t_1^2f_3^2 -
    2t_0^4t_1f_4f_1 - 2t_0^4t_1f_3f_2 - 3/2t_0^3t_1^5f_4^2 - 3t_0^3t_1^4f_4f_3 -
    3t_0^3t_1^3f_4f_2 - 3/2t_0^3t_1^3f_3^2 - 3t_0^3t_1^2f_4f_1 \\
    - 3t_0^3t_1^2f_3f_2 - t_0^3t_1f_4f_0 - t_0^3t_1f_3f_1 - 1/2t_0^3t_1f_2^2 -
    3/4t_0^2t_1^6f_4^2 - 3/2t_0^2t_1^5f_4f_3 - 3/2t_0^2t_1^4f_4f_2 -
    3/4t_0^2t_1^4f_3^2 - 3/2t_0^2t_1^3f_4f_1 \\
    - 3/2t_0^2t_1^3f_3f_2 -
    3/2t_0^2t_1^2f_4f_0 - 3/2t_0^2t_1^2f_3f_1 - 3/4t_0^2t_1^2f_2^2 + t_0t_1f_2f_0
    + 1/2t_0t_1f_1^2 + 1/2t_1f_1f_0 - 3/4f_0^2\bigg)\\
    + ab^2\bigg(-5/4t_0^7f_4^2 - 9/2t_0^6t_1f_4^2 - 2t_0^6f_4f_3 - 15/4t_0^5t_1^2f_4^2 -
    15/2t_0^5t_1f_4f_3 - 3/2t_0^5f_4f_2 - 3/4t_0^5f_3^2 - 3t_0^4t_1^3f_4^2 -
    6t_0^4t_1^2f_4f_3 - 6t_0^4t_1f_4f_2 \\
    - 3t_0^4t_1f_3^2 - t_0^4f_4f_1 -
    t_0^4f_3f_2 - 9/4t_0^3t_1^4f_4^2 - 9/2t_0^3t_1^3f_4f_3 - 9/2t_0^3t_1^2f_4f_2
    - 9/4t_0^3t_1^2f_3^2 - 9/2t_0^3t_1f_4f_1 - 9/2t_0^3t_1f_3f_2 -
    1/2t_0^3f_4f_0 \\
    - 1/2t_0^3f_3f_1 - 1/4t_0^3f_2^2 - 3/2t_0^2t_1^5f_4^2 -
    3t_0^2t_1^4f_4f_3 - 3t_0^2t_1^3f_4f_2 - 3/2t_0^2t_1^3f_3^2 -
    3t_0^2t_1^2f_4f_1 - 3t_0^2t_1^2f_3f_2 - 3t_0^2t_1f_4f_0 - 3t_0^2t_1f_3f_1 \\
    -
    3/2t_0^2t_1f_2^2 - 3/4t_0t_1^6f_4^2 - 3/2t_0t_1^5f_4f_3 - 3/2t_0t_1^4f_4f_2
    - 3/4t_0t_1^4f_3^2 - 3/2t_0t_1^3f_4f_1 - 3/2t_0t_1^3f_3f_2 -
    3/2t_0t_1^2f_4f_0 \\
    - 3/2t_0t_1^2f_3f_1 - 3/4t_0t_1^2f_2^2 - 3/2t_0t_1f_3f_0 -
    3/2t_0t_1f_2f_1 + 1/2t_0f_2f_0 + 1/4t_0f_1^2 + f_1f_0\bigg) \\
    + b^3\bigg(-7/4t_0^6f_4^2 - 3/2t_0^5t_1f_4^2 - 3t_0^5f_4f_3 - 5/4t_0^4t_1^2f_4^2 -
    5/2t_0^4t_1f_4f_3 - 5/2t_0^4f_4f_2 - 5/4t_0^4f_3^2 - t_0^3t_1^3f_4^2 -
    2t_0^3t_1^2f_4f_3 - 2t_0^3t_1f_4f_2 \\
    - t_0^3t_1f_3^2 - 2t_0^3f_4f_1 -
    2t_0^3f_3f_2 - 3/4t_0^2t_1^4f_4^2 - 3/2t_0^2t_1^3f_4f_3 -
    3/2t_0^2t_1^2f_4f_2 - 3/4t_0^2t_1^2f_3^2 - 3/2t_0^2t_1f_4f_1 -
    3/2t_0^2t_1f_3f_2 - 3/2t_0^2f_4f_0 \\
    - 3/2t_0^2f_3f_1 - 3/4t_0^2f_2^2 -
    1/2t_0t_1^5f_4^2 - t_0t_1^4f_4f_3 - t_0t_1^3f_4f_2 - 1/2t_0t_1^3f_3^2 -
    t_0t_1^2f_4f_1 - t_0t_1^2f_3f_2 - t_0t_1f_4f_0 - t_0t_1f_3f_1 \\
    - 1/2t_0t_1f_2^2 -
    t_0f_3f_0 - t_0f_2f_1 - 1/4t_1^6f_4^2 - 1/2t_1^5f_4f_3 - 1/2t_1^4f_4f_2 -
    1/4t_1^4f_3^2 - 1/2t_1^3f_4f_1 - 1/2t_1^3f_3f_2 - 1/2t_1^2f_4f_0 \\
    - 1/2t_1^2f_3f_1 - 1/4t_1^2f_2^2 - 1/2t_1f_3f_0 - 1/2t_1f_2f_1 - 1/2f_2f_0 -
    1/4f_1^2\bigg).
\end{multline*}
To satisfy (C8), we find $g_3$ as a rational function of the rest of the variables by setting the equation in (C8) to be zero and isolating for $g_3$ in the same way as earlier.
However, due to the length of the equation, it has been omitted from this paper.
The denominator has been provided below; the full definition of $g_3$ is found \thecode.
\begin{multline*}
    \textbf{denom}(g_3) = (t_0 - t_1)^3(at_0 + b)^3(at_1 + b)^2\bigg(
    (at_1 + b)(at_0 + b)^4 -
    9(at_0 + b)x_1^2 -  \\
    3(t_1 - t_0)\bigg[(at_1 + b)\bigg(f_4(t_0 + t_1)(t_0^2 + t_1^2) + f_3(t_0^2 + t_1t_0 + t_1^2) + f_2(t_0 + t_1) + f_1\bigg) - af(t_0)\bigg]\bigg).
\end{multline*}
The relationship between the free variables and dependent variables is summarized in Figure \ref{fig:DP1-vars}.
\begin{figure}[H]
    \centering
    \begin{tabular}{|l||c|c|}
        \hline
         & Free & Rational function of the others \\
        \hline
        \hline
        Coefficients of $f$: & $f_0,f_1,f_2, f_3, f_4$ & \\
        \hline
        Coefficients of $g$: & $g_4, g_5, g_6$ & $G_0, G_1, G_2, G_3$ \\
        \hline
        Coefficients of $T$: & $a, b$ & $C, D, E, H$ \\
        \hline
        Coefficients of $Q$: & $t_0$ & $X_0, Y_0$ \\
        \hline
        Coefficients of $R$: & $t_1, x_1$ & $Y_1$ \\
        \hline
    \end{tabular}
    \caption{Table showing which coefficients are free for a surface $\E\in\famone$, defined in Definition \ref{def:famone}.}
    \label{fig:DP1-vars}
\end{figure}
Since each surface contains a trisection with an additional singularity at $R$, each such trisection has genus zero from (\ref{eq:genus_formula}), meaning that each surface has Zariski dense $k$-rational points by Theorem \ref{thm:ZD}.
Comparing these equations to those found in Definition \ref{def:famone}, we see the $\E$ defined here is a member of $\famone$.
$\qed$

\begin{remark}
    Although $\famone$ has thirteen free variables, it is not a thirteen dimensional family of surfaces.
    Indeed, each surface is of the form:
    \begin{equation*}
        \E : y^2 = x^3 + f(t) x + g(t),
    \end{equation*}
    where $f,g\in k[t]$ have degree four and six, respectively.
    Without any additional constraints, this family has a total of twelve free variables, hence, any subfamily of these surfaces cannot be more than twelve-dimensional.
    Note that in Definition \ref{def:famone} the polynomial coefficients $f_0,\ldots,f_4$ and $g_4,\ldots,g_6$ are all free variables, hence, the dimension of $\famone$ is at least eight.
    The remaining free variables affect four of the coefficients of the polynomial $g$, potentially adding redundancy.
\end{remark}

\begin{example}
    The surface
    \begin{equation*}
        \E : 
        \begin{cases}
            f(t) = t^4 + t^3 + t^2 + t + 1\\
            g(t) = t^6 + t^5 + t^4 - \frac{25936}{459}t^3 + \frac{70331}{612}t^2 - \frac{201749}{2448}t +\frac{68431}{3672}
        \end{cases}
    \end{equation*}
    is a member of $\famone$ from Theorem \ref{thm:famone}, and has a trisection
    \begin{equation*}
        T : xt + x + \frac{5131}{1632}t^3 - \frac{2387}{204}t^2 + \frac{22595}{1632}t - \frac{1461}{272},
    \end{equation*}
    with a triple singularity at $Q = (4/3, 31/12, 1)$, and a double singularity at $R = (3, 29/3, 2)$.
    By equation (\ref{eq:genus_formula}), this trisection has genus zero.
\end{example}
\subsection{Unirationality and Genus Zero}
\label{sec:DP2}

Recall Definition \ref{def:fam-L} of the family $\famL$, wherein we made the assumption that $f(t_0) \neq -(at_0 + b)^4/3$.
In this section, we study when the equality holds and construct trisections on this new family of surfaces.
This allows us to prove $L$-unirationality on a subfamily of these surfaces, but at the cost of them no longer being del Pezzo surfaces of minimal degree one.
\begin{lemma}
    \label{lemma:f-t_0}
    Consider a rational elliptic surface $\E$ of the form:
    \begin{equation*}
        \E : y^2 = x^3 + f(t)x + g(t).
    \end{equation*}
    Let $Q = (x_0,y_0,t_0)$ be on $\E$, and suppose that $\E$ contains a trisection of the form:
    \begin{equation*}
        C : y = ax + bxt + ct^3 + dt^2 + et + h,
    \end{equation*}
    going through $Q$ three times.
    Further suppose that $t_0(at_0 + b)\neq 0$, and that $f(t_0) = -(at_0 + b)^4/3$.
    Then $\E$ has a minimal model which is a del Pezzo surface of degree larger than one.
\end{lemma}

\proof

Let $S,\E,C,$ and $Q$ be as in the statement of the lemma.
Recall conditions (C1-C6) from the proof of Lemma \ref{lemma:trisection}.
Looking at the conditions (C3), (C5), and (C6), using \magma{} we can write them as follows:
\begin{equation*}
    \text{(C6)} : 0 = -2\left(at_0 + b\right)^2 + 6x_0, 
\end{equation*}
which is true if and only if $x_0 = (at_0 + b)^2/3$.
Now, (C3) is equivalent to:
\begin{equation}
    \label{eq:C3-factored}
    \text{(C3)} : 0 = -2\left(at_0 + b\right)\left(ax_0t_0 + bx_0 + ct_0^3 + dt_0^2 + et_0 + h\right) + 3x_0^2 + f(t_0),
\end{equation}
however, by assumption, $f(t_0) = -(at_0 + b)^4/3$, and the condition on $x_0$ means that this is equal to $-3x_0^2$,i.e., $3x_0 + f(t_0)$ from (\ref{eq:C3-factored}) is equal to zero.
By assumption, $(at_0 + b) \neq 0$, and (C3) holds, meaning that $ax_0t_0 + bx_0 + ct_0^3 + dt_0^2 + et_0 + h = 0$.
This occurs if and only if
\begin{equation*}
    h = -ax_0t_0 - bx_0 - ct_0^3 - dt_0^2 - et_0.
\end{equation*}
Plugging these into (C5) and using \magma{}, we see that (C5) can be rewritten as:
\begin{equation*}
    \text{(C5)} : 0 = -2/3a^4t_0^3 - 2a^3bt_0^2 - 2a^2b^2t_0 - 2/3ab^3 - 6act_0^3 - 4adt_0^2 - 6bct_0^2 - 4bdt_0 - 2(at_0 + b)e + f'(t_0).
\end{equation*}
By assumption, this holds; this occurs if and only if
\begin{equation*}
    e = \frac{2/3a^4t_0^3 + 2a^3bt_0^2 + 2a^2b^2t_0 + 2/3ab^3 + 6act_0^3 + 4adt_0^2 + 6bct_0^2 + 4bdt_0 - f'(t_0)}{2(at_0 + b)}.
\end{equation*}
From here, we can plug in these conditions into the remaining conditions (C1), (C2), and (C4).
Using \magma{}, we see that the remaining conditions only depend on $a,b,t_0$, and the coefficients which define the surface $S$.
Once again, by assumption, all of these conditions are equal to zero, which can happen if and only if:
\begin{equation}
    \tag{\ref{eq:constraint-G0}}
    g_0 = 2/27(at_0 + b)^6 - t_0^6g_6 - t_0^5g_5 - t_0^4g_4 - t_0^3g_3 - t_0^2g_2 - t_0g_1,
\end{equation}
\begin{equation}
    \tag{\ref{eq:constraint-G1}}
    g_1 = -(at_0 + b)^2\left(4/3t_0^3f_4 + t_0^2f_3 + 2/3t_0f_2 + 1/3f_1\right) - 6t_0^5g_6 - 5t_0^4g_5 - 4t_0^3g_4 - 3t_0^2g_3 - 2t_0g_2,
\end{equation}
and
\begin{multline}
    \tag{\ref{eq:constraint-G2}}
    g_2 = -\frac{1}{2(at_0 + b)^2}\bigg((at_0 + b)^4\left(4t_0^2f_4 + 2t_0f_3 + 2/3f_2\right) + 2(at_0 + b)^2\left(15t_0^4g_6 + 10t_0^3g_5 + 6t_0^2g_4 + 3t_0g_3\right) \\
    - 8\left(t_0^3f_4 + 3/4t_0^2f_3 + 1/2t_0f_2 + 1/4f_1\right)^2\bigg).
\end{multline}
This shows that the corresponding code found \thecode, is a general case.
Using that code, we can see that the discriminant, $4f^3 + 27g^2$, has a factor of $(t-t_0)^3$, meaning the surface has an irreducible fiber of type III, and hence the minimal model of $\E$ will be a del Pezzo surface of degree at least two.
$\qed$

While this surface is no longer a minimal del Pezzo of degree one, we study it anyway.
It is known already that this surface is unirational, but with the condition $f(t_0) = (at_0 + b)^4/3$, we can prove it again for the following family.

\begin{definition}
    Fix a field $k$ of characteristic zero.
    Let $\famtwo$ be the family of all rational elliptic surfaces $\E$ defined by
    \begin{equation*}
        \E : y^2 = x^3 + f(t) x + g(t),
    \end{equation*}
    with polynomials $f,g\in k[t]$, $\deg(f) = 4$, $\deg(g) = 6$ defined as:
    \begin{align*}
        f(t) = F_0 + F_1t + f_2t^2 + f_3t^3 + f^4t^4, & & g(t) = G_0 + G_1t + G_2t^2 + G_3t^3 + G_4t^4 + g_5t^5 + g_6t^6.
    \end{align*}
    Here, the coefficients $f_i,g_i\in k$ are free variables, and the coefficients $F_i,G_i$ are rational functions of the aforementioned free variables along with the free variables $a,b,t_0,t_1,x_1\in k$.
    These rational functions are:
    \begin{equation}
        \label{eq:constraint-F0}
        F_0 = -\left(\frac{1}{3}(at_0 + b)^4 + t_0^4f_4 + t_0^3f_3 + t_0^2f_2 + t_0f_1\right),
    \end{equation}
    \begin{equation}
        \label{eq:constraint-F1}
        F_1 = -\frac{1}{t_0 - t_1}\left(f_4(t_0^4 - t_1^4) + \frac{1}{3}(at_0 + b)^4 + t_0^3f_3 + t_0^2f_2 - 3x_1^2 -
        t_1^3f_3 - t_1^2f_2\right),
    \end{equation}
    \begin{equation}
        \label{eq:constraint-G0}
        G_0 = \frac{2}{27}(at_0 + b)^6 - t_0^6g_6 - t_0^5g_5 - t_0^4g4 - t_0^3g_3 - t_0^2g_2 - t_0g_1,
    \end{equation}
    \begin{equation}
        \label{eq:constraint-G1}
        G_1 = -(at_0 + b)^2\left(\frac{4}{3}t_0^3f_4 + t_0^2f_3 + \frac{2}{3}t_0f_2 + \frac{1}{3}f_1\right) -
        6t_0^5g_6 - 5t_0^4g_5 - 4t_0^3g4 - 3t_0^2g_3 - 2t_0g_2,
    \end{equation}
    \begin{multline}
        \label{eq:constraint-G2}
        G_2 = -\frac{1}{2(at_0 + b)^2}\bigg((at_0 + b)^4\left(4t_0^2f_4 + 2t_0f_3 + \frac{2}{3}f_2\right) + 2(at_0 + b)^2(15t_0^4g_6 + 10t_0^3g_5 + 6t_0^2g4 + 3t_0g_3) \\
        - 8\left(t_0^3f_4 + \frac{3}{4}t_0^2f_3 + \frac{1}{2}t_0f_2 + \frac{1}{4}f_1\right)^2\bigg),
    \end{multline}
    \begin{multline}
        \label{eq:constraint-G3}
        G_3 = -\frac{1}{3(t_0 - t_1)^3(at_0 + b)^2}\bigg(\frac{1}{18}(at_0 + b)^8 - \frac{1}{3}x_1(at_0 + b)^6 + 3x_1^3(at_0 + b)^2 - \frac{9}{2}x_1^4 \\ 
        + (t_0 - t_1)^2\bigg((at_0 + b)^4\left[4f_4t_0^2 + 2f_3t_0 + \frac{2}{3}f_2\right] + 24(at_0 + b)^2(t_0 - t_1)\bigg[g_6\left(t_0^3 + \frac{3}{4}t_1t_0^2 + \frac{1}{2}t_1^2t_0 + \frac{1}{4}t_1^3\right) \\
        + g_5\left(\frac{5}{8}t_0^2 + \frac{5}{12}t_1t_0 + \frac{5}{24}t_1^2\right) + g_4\left(\frac{1}{3}t_0 + \frac{1}{6}t_1\right)\bigg] - \frac{9}{2}(t_0 - t_1)^2\bigg[f_4\left(t_0^2 + \frac{2}{3}t_0t_1 + \frac{1}{3}t_1^2\right) + 
        f_3\left(\frac{2}{3}t_0 + \frac{1}{3}t_1\right) + \frac{1}{3}f_2\bigg]^2 \\
        - x_1(at_0 + b)^2\bigg[f_4\left(t_0^2 + 2t_1t_0 + 3t_1^2\right) + f_3\left(t_0 + 2t_1\right) + f_2\bigg] - 3x_1^2\bigg[f_4\left(3t_0^2 + 2t_0t_1 + t_1^2\right) + f_3\left(2t_0 + t_1\right) + f_2\bigg]\bigg)\bigg),
    \end{multline}
    and
    \begin{multline}
        \label{eq:constraint-G4}
        G_4 = \frac{3}{(t_0 - t_1)^4(at_0 + b)^2}\bigg(\frac{1}{108}(at_0 + b)^8 - \frac{1}{9}x_1(at_0 + b)^6 + \frac{1}{2}x_1^2(at_0 + b)^4 - \frac{1}{6}(at_0 + b)^4(t_0 - t_1)^2\bigg(f_4\left(t_0^2 - 2t_0t_1 - t_1^2\right) \\
        - t_1f_3 - \frac{1}{3}f_2\bigg) + (t_0 - t_1)^4\bigg(-2(at_0 + b)^2\bigg[g_6\left(t_0^2 + t_0t_1 + \frac{1}{2}t_1^2\right) + g_5\left(\frac{1}{2}t_0 + \frac{1}{3}t_1\right)\bigg] + \frac{3}{4}\bigg[f_4\left(t_0^2 + \frac{2}{3}t_0t_1 + \frac{1}{3}t_1^2\right) \\
        + f_3\left(\frac{2}{3}t_0 + \frac{1}{3}t_1\right) + \frac{1}{3}f_2\bigg]^2\bigg) + \frac{3}{4}x_1^4 - x_1^3(at_0 + b)^2 + \frac{3}{2}x_1^2(t_0 - t_1)^2\bigg(f_4\left(t_0^2 + \frac{2}{3}t_0t_1 + \frac{1}{3}t_1^2\right) + f_3\left(\frac{2}{3}t_0 + \frac{1}{3}t_1\right) + \frac{1}{3}f_2\bigg) \\
        - \frac{1}{3}x_1(t_0 - t_1)^2(at_0 + b)^2\bigg(f_4\left(t_0^2 + 2t_0t_1 + 
        3t_1^2\right) + f_3\left(t_0 + 2t_1\right) + f_2\bigg)\bigg).
    \end{multline}

    This family is obtained by following a similar process as in Lemma \ref{lemma:f-t_0}. 
    Starting with a point $Q = (x_0,y_0,t_0)$, we construct a trisection through $Q$, with a triple singularity at $Q$.
    We then take another point, $R = (x_1,y_1,t_1)$, and force the curve to go through $R$ with multiplicity two.
    This involves three new conditions, each of which has to be zero.
    First, we set one of the two free variables in the earlier-derived two-dimensional family of trisections to a specific value; this sets one of the additional three conditions to zero.
    We then specify $f_1$ to eliminate the effect of the second free variable on the remaining two conditions.
    Afterwards, restrictions are imposed on $g_3,g_4$ so the two remaining conditions are zero, regardless of the choice of the second free variable.
    Using equation (\ref{eq:genus_formula}), we can see that each of the trisections in the now one-dimensional family will have genus zero, since they each have a triple singularity at $Q$ and now a double singularity at $R$.
    The relationship between the coefficients of the surface, the free variables not affecting the surface, and the dependent variables can be summarized in Figure \ref{fig:DP2-vars}.
\end{definition}

\begin{figure}[H]
    \centering
    \begin{tabular}{|l||c|c|c|}
        \hline
         & Free and affect $\E$ & Free and $\E$-independent & Rational function of the others \\
        \hline
        \hline
        Coefficients of $f$: & $f_2, f_3, f_4$ & & $F_0, F_1$ \\
        \hline
        Coefficients of $g$: & $g_5, g_6$ & & $G_0, G_1, G_2, G_3, G_4$ \\
        \hline
        Coefficients of $C$: & $a, b$ & $h$ & $c, d, e$ \\
        \hline
        Coefficients of $Q$: & $t_0$ & & $x_0, y_0$ \\
        \hline
        Coefficients of $R$: & $t_1, x_1$ & & $y_1$ \\
        \hline
    \end{tabular}
    \caption{Table showing which coefficients are free and which affect the surface $\E\in\famtwo$.}
    \label{fig:DP2-vars}
\end{figure}

For example, if we set all the free variables which affect the surface $\E$ to be $1$, except for $t_1$ (which must be different from $t_0$, so will be set to $2$), then we obtain the following example.

\begin{example}
    The surface $\E$ defined by
    \begin{equation*}
        \E :
        \begin{cases}
            f(t) = t^4 + t^3 + t^2 - \frac{68}{3}t + \frac{43}{3}, \\
            g(t) = t^6 + t^5 - \frac{2783}{144}t^4 + \frac{4175}{216}t^3 + \frac{501}{16}t^2 - \frac{433}{12}t + \frac{811}{108}
        \end{cases}
    \end{equation*}
    has a one dimensional family of genus zero trisection given by:
    \begin{equation*}
        T_h : y = xt + x + \left(\frac{13}{4} - \frac{1}{2}h\right)t^3 + \left(2h - \frac{103}{12}\right)t^2 + \left(\frac{8}{3} - \frac{5}{2}h\right)t + h.
    \end{equation*}
\end{example}   

\begin{theorem}[Theorem \ref{thm:intro-DP2}]
    \label{thm:DP2-unirational}
    Suppose $\E\in\famtwo$.
    Then, $\E$ has a minimal model which is a del Pezzo surface of degree two or larger.
    This minimal model is unirational, with degree of unirationality equal to three. 
\end{theorem}

\proof

The rational elliptic surface, $\E$, has a minimal model $S$ which is a del Pezzo surface of degree at least two by Lemma \ref{lemma:f-t_0}.
By construction, or by looking at the corresponding code found \thecode, we see that the trisection defined in the code passes through $Q$ with multiplicity three and $R$ with multiplicity two.
Thus, each of the trisections in the family have genus zero by (\ref{eq:genus_formula}).
This means we have a one-dimensional family of genus zero curves in $\E$, hence also in $S$, which is indexed by the free variable $h\in k$ which has no other restrictions.
This means that $S$ is unirational.
Furthermore, since each of the curves are trisections, the degree of unirationality is three.
$\qed$
\subsection{Examples}
\label{sec:examples}

As stated in Section \ref{sec:Zariski}, the difficulty of applying Theorem \ref{thm:ZD} to prove Zariski density on all elements of $\famL$ lies in finding a model of the elliptic curves $\T_R$.
To apply the theorem, one would need a Weierstrass model for infinitely many of the trisections $\T_R$, all on the same surface $\E$.
While we were not able to compute such an example yet, we were able to find families of examples, each of which has a trisection $\T_R$ for which we were able to find the Weierstrass model.

\begin{proposition}
    \label{prop:models}
    Consider the following families of rational elliptic surfaces of the form:
    \begin{equation*}
        \E : y^2 = x^3 + f(t)x + g(t).
    \end{equation*}
    Corresponding to the given values of $f$ and $g$, with parameter $\alpha$, together with two points $R$ and $Q$ written in $(x,y,t)$-coordinates.
    \begin{enumerate}[label=(\arabic*) ]
        \item $R = (13,31,7)$, $Q = (300,9000,2)$, and
        \begin{equation*}
            \hspace{\leftoffsetepsilon} \E_\alpha : 
            \begin{cases}
                f(t) = & 16886t^4 + 4t^3 + 7t^2 - 77t - 82 \\
                g(t) = & (\alpha - 421875) t^6 + (12 -4\alpha) t^5 + (4\alpha + 11)t^4 - 135t^3 + \\ 
                & (49105909849/25- 2401\alpha)t^2 + (9604\alpha - 196423631696/25)t + 196423636996/25\text{;}
            \end{cases}
            \rulefill
        \end{equation*}
        \item $R = (3 , 4 , 5)$, $Q = (4/3 , 8/3 , 2)$, and 
        \begin{equation*}
            \hspace{\leftoffsetepsilon} \E_\alpha : 
            \begin{cases}
                f(t) = & 28/3t^4 - 10t^3 - 9t^2 - 8t - 12 \\
                g(t) = & (\alpha - 1/27)t^6 + (-4\alpha + 12)t^5 + (4\alpha - 37)t^4 + 14t^3 + \\
                & (-625\alpha - 775988/243)t^2 + (2500\alpha + 3117560/243)t - (2500\alpha + 3107840/243)\text{;}
            \end{cases}
            \rulefill
        \end{equation*}
        \item $R = (1 , 1 , 2)$, $Q = (1/3 , 1/3 , 1)$, and
        \begin{equation*}
            \hspace{\leftoffsetepsilon} \E_\alpha : 
            \begin{cases}
                f(t) = & 4/3t^4 -  1\\
                g(t) = & (\alpha - 1/27)t^6 + (-2\alpha + 1)t^5 + (\alpha - 1)t^4 + \\
                & (-16\alpha - 890/27)t^2 + (32\alpha + 1753/27)t - (16\alpha + 863/27)\text{;}
            \end{cases}
            \rulefill
        \end{equation*}
        \item $R = (1 , 2 , 1)$, $Q = (4/3 , 8/3 , 2)$, and
        \begin{equation*}
            \hspace{\leftoffsetepsilon} \E_\alpha : 
            \begin{cases}
                f(t) = & 4/3t^4 - 3t^2 - 4 \\
                g(t) = & (\alpha - 1/27)t^6 + (-4\alpha + 2)t^5 + (4\alpha - 7)t^4 + 6t^3 + \\
                & (-\alpha - 8/27)t^2 + (4\alpha - 184/27)t - (4\alpha - 400/27)\text{;}  
            \end{cases}
            \rulefill
        \end{equation*}
        \item $R = (8 , 6 , 5)$, $Q = (108 , 1944 , 2)$, and
        \begin{equation*}
            \hspace{\leftoffsetepsilon} \E_\alpha : 
            \begin{cases}
                f(t) = & 2189t^4 + 2t^3 - 12t^2 + t - 2 \\
                g(t) = & (\alpha - 19683)t^6 + (-4\alpha + 2)t^5 + (4\alpha - 2)t^4 - 9t^3 + [-625\alpha + \\
                & 296597824/9]t^2 + (2500\alpha - 1186391188/9)t - (2500\alpha - 1186391440/9)\text{;}
            \end{cases}
            \rulefill
        \end{equation*}
        \item $R = (0 , 3 , 4)$, $Q = (4/3 , 8/3 , 2)$, and 
        \begin{equation*}
            \hspace{\leftoffsetepsilon} \E_\alpha : 
            \begin{cases}
                f(t) = & 19/3t^4 - 12t^3 + 5t^2 - 10t \\
                g(t) = & (\alpha - 1/27)t^6 - 4\alpha t^5 + (4\alpha + 7)t^4 - 28t^3 + \\
                & (-256\alpha - 4733/108)t^2 + (1024\alpha + 7757/27)t - (1024\alpha + 7757/27)\text{;}
            \end{cases}
            \rulefill
        \end{equation*}
        \item $R = (5 , 7 , 11)$, $Q = (12 , 72 , 3)$, and
        \begin{equation*}
            \hspace{\leftoffsetepsilon} \E_\alpha : 
            \begin{cases}
                f(t) = & 85/3t^4 - 50t^3 - 40t^2 - 38t - 39 \\
                g(t) = & (\alpha - 64/27)t^6 + (-6\alpha + 37)t^5 + (9\alpha - 191)t^4 + 176t^3 + [-14641\alpha - \\
                &  6218741/432]t^2 + (87846\alpha + 6282893/72)t - (131769\alpha + 6264101/48)\text{;}
            \end{cases}
            \rulefill
        \end{equation*}
        \item $R = (3 , 4 , 2)$, $Q = (16/3 , 64/3 , 4)$, and
        \begin{equation*}
            \hspace{\leftoffsetepsilon} \E_\alpha : 
            \begin{cases}
                f(t) = & 4/3t^4 - t^3 - 9t^2 - 11t - 4 \\
                g(t) = & (\alpha - 1/27)t^6 + (-8\alpha + 1)t^5 + (16\alpha - 4)t^4 - 12t^3 + \\ 
                & (-16\alpha + 3709/108)t^2 + (128\alpha + 1222/27)t - (256\alpha - 1012/27)\text{;}
            \end{cases}
            \rulefill
        \end{equation*}
        \item $R = (1 , 1 , 2)$, $Q = (16/3 , 64/3 , 4)$, and
        \begin{equation*}
            \hspace{\leftoffsetepsilon} \E_\alpha : 
            \begin{cases}
                f(t) = & 4/3t^4 - 4t^3 + t - 4 \\
                g(t) = & (\alpha - 1/27)t^6 + (-8\alpha + 1)t^5 + (16\alpha - 8)t^4 + 20t^3 + \\
                & (-16\alpha - 2389/54)t^2 + (128\alpha + 4372/27)t - (256\alpha + 5288/27)\text{;}
            \end{cases}
            \rulefill
        \end{equation*}
        \item $R = (4 , 1 , 5)$, $Q = (3 , 9 , 1)$, and
        \begin{equation*}
            \hspace{\leftoffsetepsilon} \E_\alpha : 
            \begin{cases}
                f(t) = & 32t^4 + t^3 - 4t^2 + 7t - 9 \\
                g(t) = & (\alpha - 27)t^6 + (-2\alpha + 3)t^5 + (\alpha - 1)t^4 + \\
                & (-625\alpha + 41613/2)t^2 + (1250\alpha - 41624)t - (625\alpha - 41631/2), \text{ and;}
            \end{cases}
            \rulefill
        \end{equation*}
        \item $R = (1 , 8 , 2)$, $Q = (196/3 , 2744/3 , 7)$, and
        \begin{equation*}
            \hspace{\leftoffsetepsilon} \E_\alpha : 
            \begin{cases}
                f(t) = & 22/3t^4 - 6t^3 - 55t^2 + t - 56 \\
                g(t) = & (\alpha - 64/27)t^6 + (-14\alpha + 5)t^5 + (49\alpha - 66)t^4 + 197t^3 + \\
                & (-16\alpha + 19423/675)t^2 + (224\alpha + 786478/675)t - (784\alpha + 1826573/675).
            \end{cases}
            \rulefill
        \end{equation*}
    \end{enumerate}
    
    All of the surfaces in each of these families contain the given points $R$ and $Q$, and contain the trisection $\T_R$ which has a triple singularity at $Q$.
    For each surface in a given family, we were able to calculate the Weierstrass model of the trisection $\T_R$, which we have omitted due to length.
    The models are available online \thecode.
    For each surface in each family, we were also able to find the coordinates of the point $P_1 + P_2$ on the original surface, which we have also omitted due to length, but are available at the same link.
    The order of the point $P_1 + P_2$ was calculated by \magma{} to be infinite.
\end{proposition}
\nocite{*}
\bibliographystyle{alpha}
\bibliography{references}
\section*{Code}

Code and examples available at: 

\url{\codelink}.

\end{document}